\documentclass[preprint]{elsarticle}
\usepackage{amssymb,bm}
\usepackage{amsmath,empheq}
\usepackage{hyperref}
\usepackage{amsthm}
\usepackage{amsbsy}
\usepackage{psfrag}
\usepackage{graphics}
\usepackage{graphicx,color}
\usepackage{tikz}
\usetikzlibrary{shapes.geometric, calc}
\usepackage{subfigure}
\usepackage{bigstrut}
\usepackage{epsfig}
\usepackage{multirow, float}
\usepackage{comment}
\usepackage{mathrsfs}
\usepackage{parskip}
\usepackage{outlines}
\usepackage{relsize}
\usepackage[normalem]{ulem}

\journal{Results in Applied Mathematics}

\theoremstyle{plain}
\newtheorem{theorem}{Theorem}[section]
\newtheorem{lemma}[theorem]{Lemma}

\theoremstyle{remark}
\newtheorem{remark}[theorem]{Remark}

\DeclareMathOperator*{\mini}{min}

\newcommand{\cQ}{\mathcal{Q}}

\newcommand{\cT}{\mathcal{T}}

\newcommand{\Ome}{\Omega}

\newcommand{\mce}{\mathcal{E}_h}
\newcommand{\mct}{\mathcal{T}_h}
\newcommand{\omct}{{\mathcal{T}}_{h,1}}
\newcommand{\twmct}{{\mathcal{T}}_{h,2}}
\newcommand{\thmct}{{\mathcal{T}}_{h,3}}

\newcommand{\bV}{\bm V}

\newcommand{\bn}{\bm n}
\newcommand{\p}{\partial}

\newcommand{\nab}{\nabla}

\newcommand{\eps}{\varepsilon}

\newcommand{\bW}{\bm W}

\newcommand{\jump}[1]{[\hspace{-0.045cm}[{#1}]\hspace{-0.045cm}]}
\newcommand{\avg}[1]{\{{#1}\}}
\newcommand{\cV}{\mathcal{V}}
\newcommand{\bcV}{{\boldsymbol{\mathcal{V}}}}
\newcommand{\bbp}{\mathbb{P}}

\definecolor{ao(english)}{rgb}{0.0, 0.5, 0.0}

\definecolor{ao(my_purple)}{rgb}{0.5, 0, 0.5}

\newcommand{\ol}{\overline}
\newcommand{\olu}{\overline{u}}
\newcommand{\oluh}{\overline{u}_h}
\newcommand{\tluh}{\tilde{u}_h}
\newcommand{\Uad}{U_{ad}}
\newcommand{\Uadh}{U_{ad,h}^{k}}
\newcommand{\dUad}{U_{ad,h}^{0}}
\newcommand{\ddUad}{U_{ad,h}^{1}}

\newcommand{\ltO}{L^{2}(\Omega)}

\newcommand{\lt}[1]{L^{2}(#1)}
\newcommand{\piho}{\Pi_h^1}
\newcommand{\pihz}{\Pi_h^0}

\newcommand{\energy}[1]{{\left\vert\kern-0.25ex\left\vert\kern-0.25ex\left\vert #1 
   \right\vert\kern-0.25ex\right\vert\kern-0.25ex\right\vert}}
\newcommand{\ds}{\displaystyle}

\newcommand{\Astar}{A^{\star}}
\newcommand{\Ahstar}{A_h^{\star}}
\definecolor{UTorange}{RGB}{247,127,0}

\begin{document}
\allowdisplaybreaks[4]
\numberwithin{equation}{section} \numberwithin{figure}{section}
\numberwithin{table}{section}

\begin{frontmatter}

\title{Convergence Analysis of a Dual-Wind Discontinuous Galerkin Method for an Elliptic Optimal Control Problem with Control Constraints}

\author[FPU]{Satyajith Bommana Boyana\corref{mycorrespondingauthor}}
\cortext[mycorrespondingauthor]{Corresponding author}
\ead{sbommanaboyana@floridapoly.edu}

\author[UNCG]{Thomas Lewis}
\ead{tllewis3@uncg.edu}

\author[WPI]{Sijing Liu}
\ead{sliu13@wpi.edu}

\author[UNCG]{Yi Zhang}
\ead{y\_zhang7@uncg.edu}

\address[FPU]{Department of Applied Mathematics, Florida Polytechnic University, Lakeland, FL, USA.}
\address[UNCG]{Department of Mathematics and Statistics, The University of North Carolina at Greensboro, Greensboro, NC, USA.}
\address[WPI]{Department of Mathematical Sciences, Worcester Polytechnic Institute, Worcester, MA, USA.}

\begin{abstract}
This paper investigates a symmetric dual-wind discontinuous Galerkin (DWDG) method for solving an elliptic optimal control problem with control constraints. The governing constraint is an elliptic partial differential equation (PDE), which is discretized using the symmetric DWDG approach. We derive error estimates in the energy norm for both the state and the adjoint state, as well as in the $L^2$ norm of the control variable. Numerical experiments are provided to demonstrate the robustness and effectiveness of the developed scheme. 
\end{abstract}


\end{frontmatter}

\section{Introduction}\label{sec:intro}
\par In this paper, we consider the following elliptic optimal control problem: 
\begin{subequations}\label{intro:sec-eocp:continuous minimization problem}
    \begin{empheq}{align}
        &\mini_{(y,u) \in H^1_0(\Omega) \times U_{ad}} &&J(y,u) := \frac{1}{2} \|y-y_d\|^2_{\lt\Ome} + \frac{\beta}{2} \|u\|^2_{\lt\Ome} \label{intro:sec-eocp:functional} \\
        &\text{subject to} &&-\Delta y = u \text{ in } \Omega, \label{intro:sec-eocp:PDE}\\
        & &&\hspace{0.31in} y = 0 \text{ on } \partial \Omega \label{intro:sec-eocp:PDE-BC}
    \end{empheq}
\end{subequations}
where $\Omega \subset \mathbb{R}^2$ is a bounded convex polygonal domain, $\beta > 0$ is a regularization parameter, $y_d \in L^2(\Omega)$ represents the desired state, and the admissible control set $U_{ad}$ is defined by 
\[
U_{ad} := \{v \in L^2(\Omega): u_a \leq v \leq u_b\}.
\]
Note that we assume $u_a < u_b$ such that $U_{ad}$ is non-empty, closed, and convex in $L^2(\Ome)$.  
This type of constraint imposed on the control variable is referred to as a box constraint. In the special case, where $u_a = -\infty$ and $u_b = \infty$, the control set reduces to $U_{ad} = L^2(\Ome)$, resulting in an optimization problem with trivial box constraints. 

\par It is well known that (see for instance, \cite{PDECO:evans2022,PDECO:renteria2021introduccion})
the state equation \eqref{intro:sec-eocp:PDE}-\eqref{intro:sec-eocp:PDE-BC} admits a unique solution $y \in H^1_0(\Omega)$ for a given $u \in L^2(\Omega)$. 
Moreover, by the continuous embedding $H^1_0(\Omega) \hookrightarrow L^2(\Omega) \hookrightarrow H^{-1}(\Omega)$, the solution operator $A: \lt\Ome \rightarrow \lt\Ome$ to  \eqref{intro:sec-eocp:PDE} - \eqref{intro:sec-eocp:PDE-BC} is linear and continuous \cite{OCP:elliptic:antil2018brief}. Therefore, for each $u$, we write the solution to \eqref{intro:sec-eocp:PDE} - \eqref{intro:sec-eocp:PDE-BC} as $y = A(u)$. Consequently, the problem \eqref{intro:sec-eocp:continuous minimization problem} reduces to 
\begin{equation} \label{intro:sec-eocp:reduced minimization problem}
    \mini_{u \in U_{ad}} J\big(A(u),u\big) = \mini_{u \in U_{ad}} \frac{1}{2} \|A(u)-y_d\|^2_{\lt\Ome} + \frac{\beta}{2} \|u\|^2_{\lt\Ome}.
\end{equation}
\par In \cite{OCP:elliptic:antil2018brief}, it was shown that \eqref{intro:sec-eocp:reduced minimization problem} admits a unique control $\ol{u} \in L^2(\Omega)$. Thus, there exists a unique state $\ol{y} = A(\ol{u}) \in H^1_0(\Omega)$ with $\left(\ol{y}, \ol{u}\right)$ uniquely satisfying \eqref{intro:sec-eocp:continuous minimization problem}.
\par Recent studies (see  \cite{OCP:elliptic:chowdhury2015framework,OCP:elliptic:casas2007using,OCP:elliptic:casas2003error,OCP:elliptic:hinze2005variational,OCP:elliptic:arada2002error} and the references therein) have extensively investigated elliptic optimal control problems (OCPs) that impose constraints on the control variable. These problems have significant applications in various engineering fields, including edge-preserving image processing \cite{OCP:elliptic:application:rudin1992nonlinear, OCP:elliptic:application:vogel2002computational}, optimizing actuator placement on piezoelectric plates to induce movement in a desired direction \cite{OCP:elliptic:application:costa2007modeling, OCP:elliptic:application:figueiredo2005piezoelectric}, and modeling total fuel consumption in vehicles \cite{OCP:elliptic:application:vossen2006l1}, among others.

\par The numerical analysis of OCPs, particularly concerning $L^2$ error estimates, has advanced significantly since the work of Falk and Geveci in the 1970s \cite{OCP:elliptic:falk1973approximation, OCP:elliptic:geveci1979approximation}, who analyzed distributed controls and Neumann boundary controls, respectively. Both authors established an $O(h)$ order of convergence using a piecewise constant approximation for the control variable. Arn{\^a}utu and Neittaanm{\"a}ki  \cite{OCP:elliptic:arnǎutu1998discretization} examined a control-constrained OCP governed by an elliptic equation in variational form within an abstract functional framework, deriving error estimates for both the optimal state and control under the assumption that \emph{a priori} error estimates for the elliptic equation hold. In \cite{OCP:elliptic:casas2003error}, Casas and Tröltzsch 
established an $O(h)$ order of convergence for the approximation of the control variable using piecewise linear, globally continuous elements in the context of linear-quadratic control problems. Later, Casas extended this result in \cite{OCP:elliptic:casas2007using} to semilinear elliptic equations and generalized objective functionals.

\par In \cite{OCP:elliptic:rosch2006error}, R{\"o}sch 
demonstrated that if both the optimal control and adjoint state are Lipschitz continuous and piecewise of class $C^2$, an improved convergence order of $O(h^{\frac32})$ could be achieved using piecewise linear approximations for the control variable in one-dimensional linear-quadratic control problems. In \cite{OCP:elliptic:hinze2005variational}, Hinze introduced a variational discretization approach and achieved  an $O(h^2)$ convergence order for the control variable. Similarly, Meyer and R{\"o}sch, in \cite{OCP:elliptic:meyer2004superconvergence}, attained the same convergence order for the control error by projecting the discrete adjoint state. R{\"o}sch and Simon derived error estimates using piecewise linear discontinuous approximations for the control variable in both 
$L^2$ and $L^{\infty}$ norms in \cite{OCP:elliptic:rosch2005linear}. More recently, Chowdhury, Gudi, and Nandakumaran \cite{OCP:elliptic:chowdhury2015framework} introduced a general framework for the error analysis of discontinuous Galerkin (DG) finite element methods applied to elliptic OCPs.

\par In this work, we propose a novel DG method based on the DG finite element differential calculus introduced in \cite{DWDG:FLN2016} to address problem \eqref{intro:sec-eocp:continuous minimization problem}. Specifically, we employ the dual-wind DG (DWDG) methods to discretize the PDE constraints given by \eqref{intro:sec-eocp:PDE} - \eqref{intro:sec-eocp:PDE-BC}. These methods have been successfully applied and analyzed in various settings, including elliptic PDEs, convection-dominated problems, as well as elliptic and parabolic obstacle problems, as evidenced in \cite{DWDG:LN2014, boyana2024convergence, EVI:DWDG:LRZ2020, PVI:DWDG:BLRZ2023}. Notably, unlike traditional DG methods, studies such as \cite{DWDG:LN2014} and \cite{PVI:DWDG:BLRZ2023} have shown that DWDG methods achieve optimal convergence rates even in the absence of a penalty term. 
To formulate a finite-dimensional problem, we define appropriate function spaces for the state variable and admissible sets for the control variable. This approach enables the formulation of the discrete Karush–Kuhn–Tucker (KKT) system and the computation of the numerical solution pair $(\ol{y}_h,\ol{u}_h)$. Given the regularity of the exact solution pair $(\ol{y},\ol{u})$, along with the discrete KKT system and the convergence analysis of DWDG methods for second-order elliptic PDEs \cite{DWDG:LN2014}, we establish convergence in the $L^2$ norm for the control variable, as well as in the energy norm for both the state and the adjoint state.

\par The structure of the paper is as follows: in Section 2, we introduce the necessary notation, review the DG finite element differential calculus framework, define various discrete operators, and discuss key properties and preliminary results that serve as the basis for later sections. In addition, we formulate the finite-dimensional optimization problem and present the corresponding discrete KKT system. Section 3 focuses on defining the energy norm, analyzing its properties, and conducting an \emph{a priori} error analysis for the control, state, and adjoint state. In Section 4, we present numerical results to validate our theoretical findings. Finally, in Section 5, we summarize the findings and discuss potential directions for future research. 
\section{Notation, the DG Calculus, and the DWDG method}\label{sec:defn}
In this section, we introduce the DG finite element differential calculus framework, establish the notation used throughout the paper, and outline key properties and results that will be useful in later sections.
\subsection{DG Operators} \label{sec:DG operators}

\subsubsection{Piecewise Sobolev Spaces and Inner Products}
\par We begin by defining the triangulation of the domain and associated sets:
\begin{itemize}
    \item Let $\mct$ denote a shape-regular simplicial triangulation of $\Ome$ \cite{FEM:BS2008,FEM:ciarlet2002finite} with mesh size $h := \max_{T \in \mct} h_T$, where $h_T$ is the diameter of the simplex $T \in \mct$.
    \item Let $\mce := \bigcup_{T \in \mct} \partial T$ denote the set of all edges in $\mct$.
    \item Let $\mce^B := \bigcup_{T \in \mct} \partial T \cap \partial \Ome$ indicate the set of boundary edges, while $\mce^I := \mce \setminus \mce^B$ represents the set of interior edges.
\end{itemize}

The set $W^{m,p}(\Ome)$ consists of all functions within $L^p(\Ome)$ whose weak derivatives up to order $m$ are also elements of $L^p(\Ome)$. In the special case where $p=2$, the space $H^m(\Ome)$, defined as $W^{m,2}(\Ome)$, becomes a Hilbert space. Furthermore, $W^{m,p}_0(\Ome)$ represents the subset of $W^{m,p}(\Ome)$ composed of functions whose traces vanish up to order $m-1$ on $\partial\Ome$. Accordingly, $H^m_0(\Ome)$ is equivalent to $W^{m,2}_0(\Ome)$. Additionally, given that $\Omega$ is a subset of $\mathbb{R}^2$, the index $i$ referenced in subsequent sections consistently assumes the values $i = 1, 2$.
\begin{itemize}
    \item Define the piecewise Sobolev spaces $W^{m.p}\left(\mathcal{T}_h\right)$ and $\bW^{m,p}(\mct)$ by:
    \begin{align*}  
        W^{m,p}(\mct)   & := \{v : v \vert_T \in W^{m,p}(T) \quad \forall T \in \mct \},                                            \\
        \bW^{m,p}(\mct) & :=\{\boldsymbol{v}: \boldsymbol{v} \vert_T \in W^{m,p}(T) \times W^{m,p}(T)  \quad \forall T \in \mct  \}.
    \end{align*}
    \item Define the inner products $\left(\cdot, \cdot\right)_{\mct}$ and $\langle \cdot, \cdot \rangle_{\mce}$, and norm $\|\cdot\|_{L^2(\mct)}$ by:
    \begin{align*}
        \ds (v, w)_{\mct} := \sum \limits_{T \in \mct}\int_T v  w \; dx,\; \; \ds \langle v, w \rangle_{\mce} := \sum \limits_{e \in \mce}\int_e v  w \; ds,\; \; \text{and} \; \; \|v\|^2_{\lt{\mct}} := (v, v)_{\mct}.
    \end{align*}
    \item Define the special subspaces $\cV_h$ and $\bcV_h$ by:
    \begin{align*}
        \cV_h:=W^{1,1}(\mct)\cap C^0(\mct) \qquad \text{and} \qquad \bcV_h:=\cV_h \times \cV_h.
    \end{align*}
\end{itemize}

\subsubsection{DG Spaces}
\begin{itemize}
    \item Define the DG space of piecewise linear polynomials $V_h$ by:
    \begin{align*}  
        V_h := \{v : v \vert_T \in \bbp_1(T) \; \; \forall T \in \mct \},
    \end{align*}
    where $\bbp_1(T)$ denotes polynomials of degree $\leq 1$ on $T$.
    \item Define the corresponding vector-valued space $\bV_h$ by:
    \begin{align*}
        \bV_h:= V_h \times V_h.
    \end{align*}
    \item[] Notice that $V_h \subset \cV_h$ and $\bV_h \subset \bcV_h$.
\end{itemize}

\subsubsection{Jump and Average Operators}
\begin{itemize}
    \item For $e = \partial T^+ \cap \partial T^- \in \mce^I$, define the jump and average operators by:
    \begin{align*}
    \jump{v}|_e:= v^+ - v^-,\qquad \avg{v}|_e:= \frac12\big( v^+ + v^-\big)\qquad \forall \, v\in \cV_h,
    \end{align*}
    where $v^\pm := v|_{T^\pm}$. Here, we denote $T^+$, $T^- \in \mct$ such that the global numbering of $T^+$ is more than that of $T^-$.
    \item For $e \in \mce^B$ (e.g., $e = \partial T^+ \cap \partial \Ome$), define the jump and average operators by:
    \begin{align*}
    \jump{v}|_e:= v^+,\qquad \avg{v}|_e:= v^+ \qquad \forall \, v\in \cV_h.
\end{align*}
\end{itemize}

\subsubsection{Trace Operators}
\begin{figure}[H]
    \centering
    \begin{minipage}{0.45\textwidth}
        \centering
        \begin{tikzpicture}[scale=1.25, transform shape]
            \draw (0,0) -- (4,0) -- (4,4) -- cycle;
            \draw (4,4) -- (0,4) -- (0,0) -- cycle;
            \draw[thick] (0,0) -- (4,4);
            \node[above,blue] at (0.25,0.5) {$e$};
            \node[above,blue] at (3.5,0.25) {$T^-$};
            \node[below,blue] at (0.5,3.75) {$T^+$};
            \draw[->,red,very thick,dashed] (2,2) -- (2,2.5);
            \node[above,red] at (2,2.4) {$\bm n_e^{(2)}$};
            \draw[->,red,very thick,dashed] (2,2) -- (1.5,2);
            \node[above right,red] at (0.75,1.7) {$\bm n_e^{(1)}$};
            \draw[->,thick] (2,2) -- (1.5,2.45);
            \node[above right,black] at (1.,2.35) {$\bm n_{e}$};
            \node[below,red] at (2.85,3.75) {$\mathcal{Q}_2^+$};
            \node[above right,red] at (3,2.5) {$\mathcal{Q}_2^-$};
            \node[below,red] at (1.25,0.95) {$\mathcal{Q}_1^+$};
            \node[above right,red] at (0.5,1.) {$\mathcal{Q}_1^-$};
        \end{tikzpicture}
    \end{minipage}
    \hfill
    \begin{minipage}{0.45\textwidth}
        \centering
        \begin{tikzpicture}[scale=1.25, transform shape]
            \draw (0,0) -- (4,0) -- (4,4) -- cycle;
            \draw (4,4) -- (0,4) -- (0,0) -- cycle;
            \draw[thick] (0,0) -- (4,4);
            \node[above,blue] at (0.25,0.5) {$e$};
            \node[above,blue] at (3.5,0.25) {$T^+$};
            \node[below,blue] at (0.5,3.75) {$T^-$};
            \draw[->,red,very thick,dashed] (2,2) -- (2,1.5);
            \node[above,red] at (2,0.9) {$\bm n_e^{(2)}$};
            \draw[->,red,very thick,dashed] (2,2) -- (2.5,2);
            \node[above right,red] at (2.45,1.7) {$\bm n_e^{(1)}$};
            \draw[->,thick] (2,2) -- (2.5,1.55);
            \node[above right,black] at (2.45,1.25) {$\bm n_{e}$};
            \node[below,red] at (2.85,3.75) {$\mathcal{Q}_2^+$};
            \node[above right,red] at (3,2.5) {$\mathcal{Q}_2^-$};
            \node[below,red] at (1.25,0.95) {$\mathcal{Q}_1^+$};
            \node[above right,red] at (0.5,1.) {$\mathcal{Q}_1^-$};
        \end{tikzpicture}
    \end{minipage}
    \caption{Trace operators $\mathcal{Q}_1^{\pm}, \mathcal{Q}_2^{\pm}$. Note that the definition is independent of the choice of $T^+$ and $T^-$.}
     \label{dgfec:sec-notation and discrete dg derivatives:fig:the trace operator}
\end{figure}
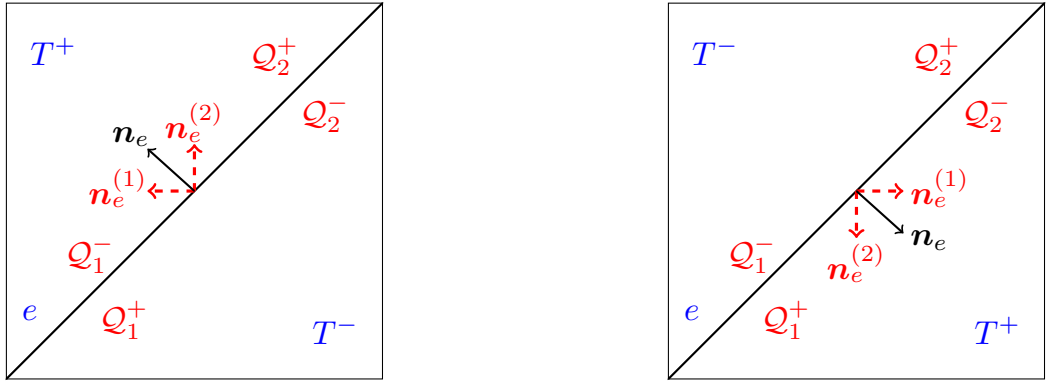

\begin{itemize}
    \item For $e \in \mce^I$, define $\bn_e = (n_e^{(1)},n_e^{(2)})^T := \bn_{T^-} \vert_{e} = -\bn_{T^+} \vert_{e}$ to be the unit normal vector. 
    \item Define the trace operators $Q_i^{\pm}(v)$ for $v \in \cV_h$ on edge $e \in \mce$ in the $x_i$ direction by:
    \begin{align*}
    \mathcal{Q}_i^+(v) :=
    \begin{cases}
        v \vert_{T^+} , & n_e^{(i)} > 0 \\
        v \vert_{T^-} , & n_e^{(i)} < 0 \\
        \{v\},          & n_e^{(i)} = 0 \\
    \end{cases}, \; \; \;
    \mathcal{Q}_i^-(v) :=
     & \begin{cases}
           v \vert_{T^-} , & n_e^{(i)} > 0 \\
           v \vert_{T^+} , & n_e^{(i)} < 0 \\
           \{v\},          & n_e^{(i)} = 0 \\
       \end{cases}.
    \end{align*}
    This definition allows us to interpret $Q_i^+(v)$ and $Q_i^-(v)$ (see Figure \ref{dgfec:sec-notation and discrete dg derivatives:fig:the trace operator}) as ``forward'' and ``backward'' limits in the $x_i$ direction on $e \in \mce^I$.
    \item For $e = \partial T^+ \cap \partial \Ome \in \mce^B$, we define $Q_i^+(v)$ by:
    \begin{align*}
        Q_i^{\pm}(v) := v^+.
    \end{align*}
\end{itemize}

\subsubsection{Discrete Partial Derivatives and Gradient Operators}
\par With the trace operators defined above, we introduce the discrete partial derivative $\partial_{h,x_i}^{\pm}:\cV_h \rightarrow V_h$ for any $v \in \cV_h$ as 

\begin{align}
        \bigl(\partial_{h,x_i}^\pm v,\varphi_h\bigr)_{\cT_h}
        &:= \bigl\langle \cQ_i^\pm (v) n^{(i)}, \jump{\varphi_h} \bigr\rangle_{\mce}
        -\bigl(v, \partial_{x_i}\varphi_h \bigr)_{\cT_h} \qquad \forall \varphi_h \in V_h. \label{dgfec:sec-notation and discrete dg derivatives:DG derivative w/o bc}  
\end{align}

\par 
Then we can naturally define the discrete gradient operator as follows. For any $v \in \cV_h$, we define
\begin{align*}
    \nab_h^\pm v := (\p_{h,x_1}^\pm v, \p_{h,x_2}^\pm v). 
\end{align*}
\subsection{The Discrete Problem} \label{sec:optimality condn}
We first present the first-order optimality conditions for the continuous problem \eqref{intro:sec-eocp:continuous minimization problem}/\eqref{intro:sec-eocp:reduced minimization problem}. It is standard that \eqref{intro:sec-eocp:continuous minimization problem}/\eqref{intro:sec-eocp:reduced minimization problem} is equivalent to the following  variational inequality:
\begin{align}\label{eocp:sec-optimality conditions:continuous VI}
    \big(A^{\star}\big(A(\ol{u})-y_d\big)+\beta \ol{u},u-\ol{u}\big)_{\lt\Ome} \; \;  \geq  \hspace{0.1in} 0 \qquad \forall u \in \Uad,
\end{align}
where $A^{\star}: \lt\Ome \rightarrow \lt\Ome$ is the adjoint operator of A.  Since for every $u \in \lt\Ome$ there exists a unique $y = A(u) \in \lt\Ome$, we define $p:=A^{\star}\big(A(u)-y_d\big) = A^{\star}(y - y_d) .$ As a result, we have the \emph{adjoint equation},
\begin{subequations}\label{eocp:sec-optimality conditions:adj eqn}
    \begin{empheq}{align}
        &-\Delta p = y - y_d \quad \text{ in } \Omega, \label{theADJ}\\
        &\hspace{0.31in} p = 0 \quad \hspace{0.275in} \text{ on } \Gamma.
    \end{empheq}
\end{subequations}
\par The solution to \eqref{eocp:sec-optimality conditions:adj eqn} is the \emph{adjoint state} $p \in H^1_0(\Ome).$ 
Finally, we have the following first-order optimality conditions for the solution $(\bar{y}, \bar{u}, \bar{p}) \in H^1_0(\Omega) \times H^1_0(\Omega) \times U_{ad}$: 
\begin{subequations}\label{eocp:sec-optimality conditions:continuous KKT}
    \begin{empheq}{align}
        (\nab\ol{y} , \nab{v})_{\lt{\Ome}} -
        (\ol{u} , v)_{\lt{\Ome}} &= 0 \qquad \forall v \in H^1_0(\Ome) \label{eocp:sec-optimality conditions:continuous KKT PDE weakform},  \\
        (\nab\ol{p} , \nab{v})_{\lt{\Ome}} - (\ol{y} - y_d, v)_{\lt{\Ome}} &= 0 \qquad \forall v \in H^1_0(\Ome) \label{eocp:sec-optimality conditions:continuous KKT adj weakform}, \\
        \big(\ol{p}+\beta \ol{u},u-\ol{u}\big)_{\lt\Ome} \; \; &\geq 0 \qquad \forall u \in \Uad \label{eocp:sec-optimality conditions:continuous KKT VI}.
    \end{empheq}
\end{subequations}
\par The coupled system \eqref{eocp:sec-optimality conditions:continuous KKT} is the first order necessary and sufficient optimality system for solving \eqref{intro:sec-eocp:continuous minimization problem}/\eqref{intro:sec-eocp:reduced minimization problem} because $J: H^1_0(\Ome) \times U_{ad} \rightarrow \mathbb{R}$ is a convex functional.  Furthermore, in view of \eqref{eocp:sec-optimality conditions:continuous KKT PDE weakform}-\eqref{eocp:sec-optimality conditions:continuous KKT adj weakform} and the convexity of the domain, we can guarantee the regularity of $\ol{y} \in H^2(\Ome)$ and $\ol{p} \in H^2(\Ome)$ (see \cite{OCP:elliptic:savare1998regularity, OCP:elliptic:troltzsch2010finite}). Also, using \eqref{eocp:sec-optimality conditions:continuous KKT VI}, it can be shown that $\ol{u} \in H^1(\Ome)$ (cf. \cite{OCP:elliptic:casas2003error} and the references therein).

Now we define the discrete admissible sets for the control variable by
\begin{align*}
    U_{ad,h}^{k} := \{v \in U_{ad}: v|_{T} \in \mathbb{P}_k(T) \; \; \forall T \in \mct \} \qquad \text{for}\; k \in \{0,1\}.
\end{align*} 
Then the discrete problem is as follows:
\begin{subequations}\label{eocp:sec-optimality conditions:discrete minimization}
    \begin{empheq}{align}
        &\mini_{(y_h,u_h) \in V_h \times U_{ad,h}^k} &&J_h(y_h,u_h) := \frac{1}{2} \|y_h-y_{d}\|^2 + \frac{\beta}{2} \|u_h\|^2
        \label{eocp:sec-optimality conditions:discrete functional} \\
        &\text{subject to } &&a_h(y_h,v_h) = (u_h,v_h)_{\cT_h} \quad \forall v_h \in V_h
        \label{eocp:sec-optimality conditions:discrete weak form}
    \end{empheq}
\end{subequations}
where
\begin{align}\label{eq:ah}
     & a_h(v_h,w_h) :=\ds \frac{1}{2} \Big( \big(\nabla_{h,0}^+ v_h, \nabla_{h,0}^+ w_h\big)_{\cT_h} + \big(\nabla_{h,0}^- v_h, \nabla_{h,0}^- w_h\big)_{\cT_h} \Big) + \bigg\langle  \frac{\gamma_e}{h_e} \jump{v_h},\jump{w_h} \bigg\rangle_{\mce}
\end{align}
for $\gamma_e$, a parameter defined on $e \in \mathcal{E}_h$ that will be determined later. 

\par Similar to the continuous case, there exists a unique solution pair: $(\ol{y}_h,\ol{u}_h) \in V_h \times U_{ad,h}^k$ satisfying \eqref{eocp:sec-optimality conditions:discrete minimization}. 
We define a discrete solution operator $A_h: U^k_{ad,h} \rightarrow V_h$ to \eqref{eocp:sec-optimality conditions:discrete weak form}   
and the discrete adjoint state $\ol{p}_h:=A_h^{\star}\left(A_h(\ol{u}_h) - y_{d}\right) = A_h^{\star}\left(\ol{y}_h - y_{d}\right) \in V_h$, where $A^{\star}_h$ denotes the adjoint operator of $A_h$. 
Finally, we have a discrete coupled system satisfied by $(\ol{y}_h, \ol{u}_h, \ol{p}_h) \in V_h \times U^k_{ad,h} \times V_h$:
\begin{subequations}\label{eocp:sec-optimality conditions:discrete KKT}
    \begin{empheq}{align}
        a_h(\ol{y}_h,v_h) &= (\ol{u}_h,v_h)_{\cT_h} \; \; &&\forall v_h \in V_h \label{eocp:sec-optimality conditions:discrete KKT PDE weakform},  \\
        a_h(\ol{p}_h,v_h) &= (\ol{y}_h - y_{d},v_h)_{\cT_h} \; \;  &&\forall v_h \in V_h\label{eocp:sec-optimality conditions:discrete KKT adj weakform}, \\
        \big(\ol{p}_h+\beta \ol{u}_h,u_h-\ol{u}_h\big)_{\cT_h} &\geq 0 \qquad &&\forall u_h \in \Uadh  \label{eocp:sec-optimality conditions:discrete KKT VI}.
    \end{empheq}
\end{subequations}
Similar to \eqref{eocp:sec-optimality conditions:continuous VI}, we have the following discrete variational inequality:
\begin{align}
        \big(A_h^{\star}(A_h\ol{u}_h-y_{d})+\beta \ol{u}_h,u_h-\ol{u}_h\big)_{\lt\Ome}  &\geq 0   \qquad \forall u_h \in U_{ad,h}^{k}. \label{sec:err estimate:discrete VI operator form}
\end{align}
\section{An \emph{a Priori} Error Estimate of the Control Variable}\label{sec:a priori analysis}
In this section, we provide an \emph{a priori} error estimate for the control variable. We follow the standard approach in \cite{OCP:elliptic:casas2003error}. The key is to construct suitable projection operators. We also need the error estimates of DWDG methods for Poisson equations (cf. \cite{DWDG:LN2014}).

\subsection{Preliminary Estimates}

We introduce the following notations:
\begin{subequations}
    \begin{empheq}{align}
        \|v_h\|^2_{1,h} &:= \frac12 \Big( \|\nab_{h,0}^+ v_h\|_{L^2(\Ome)}^2
        + \|\nab_{h,0}^- v_h\|_{L^2(\Ome)}^2\Big) \qquad \forall \, v_h \in V_h, \\
        \energy{v_h} &:= \|v_h\|^2_{1,h} + \sum_{e\in \mce} \frac{\gamma_e}{h_e} \big\|\jump{v_h}\big\|_{L^2(e)}^2 \qquad \forall \, v_h \in V_h. \label{sec:a priori analysis:energynorm}
    \end{empheq}
\end{subequations}

\begin{theorem}\label{sec:a priori analysis:thm:estimate controlling the jump}
    Let $\gamma_{\min}:= \min_{e\in \mce} \gamma_e$. Then
    \begin{align}
        \gamma_{\min}\sum_{e\in \mce} h_e^{-1}\big\|\jump{v_h}\big\|_{L^2(e)}^2 \le \energy{v_h}^2 \qquad \forall \, v_h \in V_h. \label{sec:a priori analysis:controlling the jump by energy norm pos penalty}
    \end{align}
    provided $\gamma_{\min} > 0$. Moreover, if the triangulation $\mct$ is quasi-uniform and each $T \in \mct$ has at most one boundary edge, then there exists a constant $C_* > 0$ independent of h and $\gamma_e \hspace{0.05in} \forall e \in \mce$ such that
    \begin{align}
        \big(C_*+\gamma_{\min}\big)\sum_{e\in \mce} h_e^{-1}\big\|\jump{v_h}\big\|_{L^2(e)}^2 \le \energy{v_h}^2 \qquad \forall \, v_h\in V_h.
            \label{sec:a priori analysis:controlling the jump by energy norm neg penalty} 
    \end{align}
\end{theorem}
\begin{proof}
    The proof of Theorem \ref{sec:a priori analysis:thm:estimate controlling the jump} can be found in \cite{DWDG:LN2014}.
\end{proof} 

Next, we note that the following relationship holds between the classical gradient and the DG discrete gradient. The proof is provided in \cite[Lemma 4.1]{EVI:DWDG:LRZ2020}.
\begin{lemma}\label{sec:a priori analysis:lemma:regulargrad_DGgrad}
    For $\gamma_{\min} > 0$, we have
    \begin{align}
         \|\nab v_h\|_{L^2(\mct)}^2\le C\left(1+ \frac{1}{\gamma_{\min}}\right)\energy{v_h}^2 \qquad \forall v_h\in V_h. \label{sec:a priori analysis:controlling regular grad by DG grad pos}
    \end{align}
    Further, if $-C_* < \gamma_{\min} \le 0$ and the triangulation $\mct$ is quasi-uniform and each simplex in the triangulation has at most one boundary edge, then
    \begin{align}
        \|\nab v_h\|_{L^2(\mct)}^2 \le C\left(1 +\frac{1+|\gamma_{\min}|}{C_*+ \gamma_{\min}}\right)\energy{v_h}^2 \qquad \forall v_h\in V_h. \label{sec:a priori analysis:controlling regular grad by DG grad neg}
    \end{align}
\end{lemma}
We then have the following discrete Poincar{\' e} inequality \cite{PVI:DWDG:BLRZ2023}.
\begin{lemma}\label{sec:a priori analysis:lemma:DWDG Poincare}
  There exists a positive constant $C$ independent of $h$ such that
  \begin{align}
      \| v_h \|^2_{L^2(\Omega)} \leq C \energy{v_h}^2 \qquad \forall v_h \, \in V_h.
  \end{align}
\end{lemma}

\subsection{Estimates on $A_h$ and $A_h^{\star}$}\label{sec:aas}
For any $v\in L^2(\Omega)$, it is easy to see that $A_hv$ is the DWDG approximation of the variable $Av\in H^1_0(\Omega)$, which satisfies a Poisson equation on the convex domain. Then, we immediately have the following estimate from \cite{DWDG:LN2014}:
\begin{equation}\label{eq:Ahl2}
    \|Av-A_hv\|_{L^2(\Omega)}\le Ch^2\|v\|_{L^2(\Omega)}.
\end{equation}
It follows from the Poincar{\' e} inequality that, for any $v\in L^2(\Omega)$,
\begin{equation}
    \|Av\|_{L^2(\Omega)}\le C \|\nabla Av\|_{L^2(\cT_h)}\le C\|v\|_{L^2(\Omega)}.
\end{equation}
We also have, for $A_hv\in V_h$,
\begin{equation}
\begin{aligned}
    \|A_hv\|^2_{L^2(\Omega)}&\le C\energy{A_hv}^2=C a_h(A_hv,A_hv)=C(v,A_hv)_{L^2(\Omega)}\le C \|v\|_{L^2(\Omega)}\|A_hv\|_{L^2(\Omega)}
\end{aligned}
\end{equation}
by Lemma \ref{sec:a priori analysis:lemma:DWDG Poincare}, \eqref{eq:ah}, \eqref{sec:a priori analysis:energynorm} and \eqref{eocp:sec-optimality conditions:discrete weak form}. We then obtain
\begin{equation}
\begin{aligned}
    \|A_hv\|_{L^2(\Omega)}\le C\|v\|_{L^2(\Omega)}.
\end{aligned}
\end{equation}
Similarly, $A^{\star}$ and $A_h^{\star}$ represent the solution operators of the dual problem of \eqref{intro:sec-eocp:PDE} and \eqref{eocp:sec-optimality conditions:discrete weak form}, respectively. We can get the following for any $v\in L^2(\Omega)$:
\begin{alignat}{3}
    \|A^{\star}v-A^{\star}_hv\|_{L^2(\Omega)}&\le Ch^2\|v\|_{L^2(\Omega)},\\
    \|A^{\star}v\|_{L^2(\Omega)}&\le C\|v\|_{L^2(\Omega)},\\
    \|A^{\star}_hv\|_{L^2(\Omega)}&\le C\|v\|_{L^2(\Omega)}.
\end{alignat}

\begin{remark}
    The operators $A$ and $A^\star$ (resp., $A_h$ and $A_h^\star$) are identical in this work since our PDE constraint is a symmetric problem. However, we use different notations to distinguish them, allowing us to track the different roles these operators play. Moreover, this distinction makes it straightforward to extend our theory to non-symmetric PDE constraints.
\end{remark}

\subsection{$\mathbb{P}_0$ Approximation of the Control}\label{sec:p0 approx control}

In this section, we provide an error estimate on the control variable in the $L^2$ norm when the finite-dimensional admissible set is $U_{ad,h}^{0}$. 
Define $\pihz: \Uad\longrightarrow \dUad$ such that 
    \begin{align*}
        \pihz(v)|_{T} = \int_{T}\frac{v}{\text{meas}(T)} \; dx \qquad \forall T \in \mathcal{T}_h.
    \end{align*}
The operator $\pihz$ is an $L^2$ projection of $\Uad$ onto $\dUad$,
and we have the following standard estimate \cite{FEM:ciarlet1991handbook,FEM:BS2008}:
\begin{align} 
    \|\ol{u} - \pihz(\ol{u})\|_{L^2(\Omega)} \leq Ch \|\nab \ol{u}\|_{L^2(\Omega)} \label{sec:p0 approx control:L2 proj estimate}.
\end{align}

\begin{theorem}\label{sec:p0 approx control:thm:L2 err in exact and discrete control}
    Let $\ol{u} \in U_{ad} \cap H^1(\Omega)$ and $\ol{u}_h \in U_{ad,h}^0 \subset U_{ad}$ be the solutions of the problems - \eqref{intro:sec-eocp:continuous minimization problem} and \eqref{eocp:sec-optimality conditions:discrete minimization}, respectively. Then, there exists a constant C that depends on $\|\ol{u}\|_{H^1(\Omega)}$ and $\|y_d\|_{L^2(\Omega)}$ and is independent of $h$ such that \\
    \begin{align*}
        \|\ol{u} - \ol{u}_h \|_{L^2(\Omega)} \leq Ch.
    \end{align*}
\end{theorem}
\begin{proof}
    \par By considering the variational inequality \eqref{eocp:sec-optimality conditions:continuous KKT VI} and from Subsection \ref{sec:optimality condn}, we have
    \begin{align}
        \left(A^{\star}(A\ol{u}-y_d)+\beta \ol{u},u-\ol{u}\right)_{\lt\Ome}  &\geq 0 \quad \forall u \in \Uad \label{sec:p0 approx control:continuous VI}.
\   \end{align}
    \par Likewise, from \eqref{sec:err estimate:discrete VI operator form}, we have
    \begin{align}
            \left(A_h^{\star}(A_h\ol{u}_h-y_d)+\beta \ol{u}_h,u_h-\ol{u}_h\right)_{\lt\Ome}  &\geq 0   \quad \forall u_h \in U_{ad,h}^{0} \label{sec:err estimate:discrete VI operator form-p0}.
    \end{align}
    \par Since $\ol{u}_h, \pihz(\ol{u}) \in \dUad \subset U_{ad}$, upon replacing $u$ by $\ol{u}_h$ and $u_h$ by $\pihz(\ol{u})$ in \eqref{sec:p0 approx control:continuous VI} and \eqref{sec:err estimate:discrete VI operator form-p0}, respectively, we have
    \begin{subequations}
        \begin{align}
            0 &\leq \big(A^{\star}(A\ol{u}-y_d)+\beta \ol{u},\ol{u}_h-\ol{u}\big)_{\lt\Ome} \label{sec:p0 approx control:p0 proof continuous VI}, \\
            0 &\leq
            \big(A_h^{\star}(A_h\ol{u}_h-y_d)+\beta \ol{u}_h,\pihz(\ol{u})-\ol{u}_h\big)_{\lt\Ome} \notag  \\ 
            &= 
            \big(A_h^{\star}(A_h\ol{u}_h-y_d)+\beta \ol{u}_h,\pihz(\ol{u})-\ol{u}+\ol{u}-\ol{u}_h\big)_{\lt\Ome} \notag \label{sec:p0 approx control:p0 proof discrete VI}\\
            &=
            \big(A_h^{\star}(A_h\ol{u}_h-y_d)+\beta \ol{u}_h,\pihz(\ol{u})-\ol{u}\big)_{\lt\Ome} \notag \\
            &\hspace{0.25in}- 
            \big(A_h^{\star}(A_h\ol{u}_h-y_d)+\beta \ol{u}_h,\ol{u}_h - \ol{u}\big)_{\lt\Ome}.
        \end{align}
    \end{subequations}
    \par Upon adding \eqref{sec:p0 approx control:p0 proof continuous VI} and \eqref{sec:p0 approx control:p0 proof discrete VI}, we obtain
    \begin{align*}
        0 &\leq \big(A^{\star}\left(A\ol{u}-y_d\right)-A_h^{\star}\left(A_h\ol{u}_h-y_d\right), \ol{u}_h-\ol{u}\big)_{\lt{\Ome}} \\
        &\hspace{0.175in}+\big(A_h^{\star}\left(A_h\ol{u}_h-y_d\right),\pihz(\ol{u})-\ol{u} \big)_{\lt{\Ome}} +\big(\beta \ol{u}_h,\pihz(\ol{u})-\ol{u} \big)_{\lt{\Ome}}\\
        &\hspace{0.175in}-\beta \left(\ol{u}_h-\ol{u} ,\ol{u}_h-\ol{u}\right)_{\lt\Ome}.
    \end{align*}
    \par By reordering and noting that $\big(\beta \ol{u}_h, \pihz(\ol{u}) - \ol{u}\big)_{\lt\Ome} = 0$ due to orthogonality, we arrive at
    \begin{align}
        \beta \big(\ol{u}_h-\ol{u} ,\ol{u}_h-\ol{u}\big)_{\lt\Ome} 
         &\leq 
         \big(A^{\star}\left(A\ol{u}-y_d\right)-A_h^{\star}\left(A_h\ol{u}_h-y_d\right), \ol{u}_h-\ol{u}\big)_{\lt{\Ome}} + \big(A_h^{\star}\left(A_h\ol{u}_h-y_d\right),\pihz(\ol{u})-\ol{u} \big)_{\lt{\Ome}} \notag \\
        &= \underbrace{\big(A^{\star}\left(A\ol{u}\right)-A_h^{\star}\left(A_h\ol{u}_h\right), \ol{u}_h-\ol{u}\big)_{\lt{\Ome}}}_{T_1} \notag \\
        &\hspace{0.2in}+ \underbrace{\left(\Ahstar y_d-\Astar y_d,  \ol{u}_h-\ol{u}\right)_{\lt{\Ome}}}_{T_2} \notag \\
        &\hspace{0.2in}+ \underbrace{\big(A^\star_h\left(A_h\ol{u}_h - A_h\ol{u}\right),\pihz(\ol{u})-\ol{u}\big)_{\ltO}}_{T_3} \notag \\
        &\hspace{0.2in}+\underbrace{\big(A^\star_h\left(A_h\ol{u} - A\ol{u}\right),\pihz(\ol{u})-\ol{u}\big)_{\ltO}}_{T_4} \notag \\
        &\hspace{0.2in}+ \underbrace{\big(A^\star_h\left(A\ol{u} - y_d\right)-A^\star\left(A\ol{u} - y_d\right),\pihz(\ol{u})-\ol{u}\big)_{\ltO}}_{T_5} \notag \\
        &\hspace{0.2in}+ \underbrace{\big(A^\star\left(A\ol{u} - y_d\right),\pihz(\ol{u})-\ol{u}\big)_{\ltO}}_{T_6}. \label{ineq:T1-6}
    \end{align}
We now estimate $T_1-T_6$ term by term, where we repeatedly use the estimates established in Section \ref{sec:aas}.
        \begin{alignat}{2}
            &T_1 & & = \big(A^{\star}\left(A\ol{u}\right)-A_h^{\star}\left(A_h\ol{u}_h\right), \ol{u}_h-\ol{u}\big)_{\lt{\Ome}} \notag \\
            &&&=\big(A^{\star}\left(A\ol{u}\right)-A_h^{\star}\left(A_h\ol{u}_h\right), \ol{u}_h-\ol{u}\big)_{\lt{\Ome}} 
            - \big(A_h^{\star}\left(A_h\ol{u}\right),\ol{u}_h-\ol{u}\big)_{\ltO} 
            + \big(A_h^{\star}\left(A_h\ol{u}\right),\ol{u}_h-\ol{u}\big)_{\ltO} \notag \\
            &&&= \big(\left(A^\star A - A_h^\star A_h\right)(\ol{u}),\ol{u}_h - \ol{u}\big)_{\ltO}
                    +\big(A_h^{\star}\left(A_h\left(\ol{u}-\ol{u}_h\right)\right), \ol{u}_h-\ol{u}\big)_{\lt{\Ome}} \notag \\
            &&&= \big(\left(A^\star A - A_h^\star A_h\right)(\ol{u}),\ol{u}_h - \ol{u}\big)_{\ltO}
                    -\big(A_h\left(\ol{u}-\ol{u}_h\right), A_h\left(\ol{u}-\ol{u}_h\right)\big)_{\ltO} \notag \\
            &&&= \big(\left(A^\star A - A_h^\star A_h\right)(\ol{u}),\ol{u}_h - \ol{u}\big)_{\ltO}
                    -\big\|A_h\left(\ol{u}-\ol{u}_h\right)\big\|^2_{\ltO} \notag \\
            &&&\leq \big(\left(A^\star A - A_h^\star A_h\right)(\ol{u}),\ol{u}_h - \ol{u}\big)_{\ltO} \notag \\
            &&&= \big(\left(A^\star A - A_h^\star A_h\right)(\ol{u}),\ol{u}_h - \ol{u}\big)_{\ltO} 
                    +\big(\left(A_h^\star A\right)(\ol{u}),\ol{u}_h - \ol{u}\big)_{\ltO}
                    -\big(\left(A_h^\star A\right)(\ol{u}),\ol{u}_h - \ol{u}\big)_{\ltO} \notag \\
            &&&= \big(\left(A^\star-A^\star_h\right)(A\ol{u}),\ol{u}_h - \ol{u}\big)_{\ltO}
                    +\big(A^\star_h\left(A\ol{u}-A_h\ol{u}\right),\ol{u}_h - \ol{u}\big)_{\ltO} \notag \\
            &&&\leq \dfrac{C}{\eps_1}\big\|\left(A^\star-A^\star_h\right)(A\ol{u})\big\|^2_{\ltO} + C\eps_1\big\|\ol{u}_h - \ol{u}\big\|^2_{\ltO}
            + \dfrac{C}{\eps_2}\big\|\left(A-A_h\right)\ol{u}\big\|^2_{\ltO} + C\eps_2\big\|\ol{u}_h - \ol{u}\big\|^2_{\ltO} \notag \\
            &&&\leq Ch^4 \big\|\ol{u}\big\|_{\ltO}^2 + Ch^4\big\|\ol{u}\big\|_{\ltO}^2 + C(\eps_1+\eps_2) \big\|\ol{u}_h - \ol{u}\big\|^2_{\ltO}. \label{boundT1} \\ \notag \\
            &T_2 & & = \big(\Ahstar y_d-\Astar y_d,  \ol{u}_h-\ol{u}\big)_{\lt{\Ome}} \notag \\
            &&&=\big(\left(A^\star - A^\star_h\right)y_d,\ol{u}_h - \ol{u}\big)_{\ltO}  \notag\\
            &&&\leq \dfrac{C}{\eps_3} \big\|\left(A^\star - A^\star_h\right)y_d\big\|^2_{\ltO} + C \eps_3\big\|\ol{u}_h - \ol{u}\big\|^2_{\ltO} \notag\\
            &&&\leq Ch^4 \big\|y_d\big\|_{\ltO}^2 + C\eps_3 \big\|\ol{u}_h - \ol{u}\big\|^2_{\ltO}. \label{boundT2} \\ \notag\\
            &T_3 & & = \big(A^\star_h\left(A_h\ol{u}_h - A_h\ol{u}\right),\pihz(\ol{u})-\ol{u}\big)_{\ltO} \notag\\
            &&&=\big(\left(A^\star_h A_h\right) (\ol{u}_h - \ol{u}),\pihz(\ol{u})-\ol{u}\big)_{\ltO}  \notag\\
            &&&\leq C\eps_4\big\|A^\star_h A_h (\ol{u}_h - \ol{u})\big\|^2_{\ltO} + \dfrac{C}{\eps_4} \big\|\pihz(\ol{u})-\ol{u}\big\|^2_{\ltO} \notag\\
            &&&\leq C\eps_4 \big\|\ol{u}_h - \ol{u}\big\|^2_{\ltO} + \dfrac{C}{\eps_4} h^2 \big\|\ol{u}\big\|^2_{H^1(\Ome)}. \label{boundT3} \\ \notag\\
            &T_4 & & =\big(A^\star_h\left(A_h\ol{u} - A\ol{u}\right),\pihz(\ol{u})-\ol{u}\big)_{\ltO} \notag \\
            &&&\leq \|A^\star_h\left(A_h\ol{u} - A\ol{u}\right)\|_{\ltO} \; \|\pihz(\ol{u})-\ol{u}\|_{\ltO}  \notag \\
            &&&\leq C\|\left(A-A_h\right)\ol{u}\|_{\ltO} \; \|\pihz(\ol{u})-\ol{u}\|_{\ltO} \notag \\
            &&&\leq Ch^3 \|\ol{u}\|^2_{H^1(\Ome)}. \label{boundT4} \\ \notag\\
            &T_5 & & = \big(A^\star_h\left(A\ol{u} - y_d\right)-A^\star\left(A\ol{u} - y_d\right),\pihz(\ol{u})-\ol{u}\big)_{\ltO} \notag\\
            &&&= \left(A^\star_h\left(\ol{y} - y_d\right)-A^\star\left(\ol{y} - y_d\right),\pihz(\ol{u})-\ol{u}\right)_{\ltO} \notag \\
            &&&= \big(\left(A^\star_h - A^\star\right)\ol{y}, \pihz(\ol{u})-\ol{u}\big)_{\ltO} + \big(\left(A^\star - A^\star_h\right)y_d, \pihz(\ol{u})-\ol{u}\big)_{\ltO} \notag \\
             &&&\leq \|\left(A^\star - A^\star_h\right)\ol{y}\|_{\ltO} \; \|\pihz(\ol{u})-\ol{u}\|_{\ltO} + \|\left(A^\star - A^\star_h\right)y_d\|_{\ltO} \; \|\pihz(\ol{u})-\ol{u}\|_{\ltO} \notag \\
             &&&\leq Ch^3\|\ol{u}\|_{\ltO} \; \|\ol{u}\|_{H^1(\Ome)} + Ch^3\|y_d\|_{\ltO} \; \|\ol{u}\|_{H^1(\Ome)}. \label{boundT5} \\ \notag \\
             & \text{Before } & &  \text{we proceed to estimate $T_6$, notice that } \pihz(\ol{p}) \in V_h \text{ and }\left(\pihz(\ol{p}),\pihz(\ol{u})-\ol{u}\right)_{\ltO} = 0. \text{ Therefore}  \notag \\
            &T_6 & & = \big(A^\star\left(A\ol{u} - y_d\right),\pihz(\ol{u})-\ol{u}\big)_{\ltO} \notag \\
             &&& = \big(\ol{p} - \pihz(\ol{p}),\pihz(\ol{u})-\ol{u}\big)_{\ltO} 
            + \big(\pihz(\ol{p}),\pihz(\ol{u})-\ol{u}\big)_{\ltO} \notag \\
            &&&\leq \|\ol{p} - \pihz(\ol{p})\|_{\ltO} \; \|\pihz(\ol{u})-\ol{u}\|_{\ltO} \notag \\
            &&&\leq Ch^2 (\|\ol{u}\|_{L^2(\Ome)} + \|y_d\|_{L^2(\Ome)}) \; \|\ol{u}\|_{H^1(\Ome)}. \label{boundT6}
        \end{alignat}
        Finally, by collecting the estimates \eqref{boundT1}-\eqref{boundT6}, by appropriately choosing $\eps_1, \eps_2, \eps_3, \eps_4$ to be sufficiently small and possibly dependant on $\beta$, and by using the inequality \eqref{ineq:T1-6}, we can conclude
        \begin{align*}
        \|\ol{u} - \ol{u}_h \|_{L^2(\Omega)} \leq C (\|\ol{u}\|_{H^1(\Ome)} + \|y_d\|_{L^2(\Ome)}) h.
    \end{align*}


\end{proof}
\subsection{$\mathbb{P}_1$ Approximation of the Control}\label{sec:p1 approx control}

\par It was shown in \cite{PVI:kinderlehrer2000introduction} that the optimal control $\bar{u} \in W^{1, \infty}(\Omega)$. In this subsection, an improved bound on the error associated with the control variable is established in the $L^2$ norm when the finite-dimensional admissible set is $U_{ad,h}^{1}.$ 

\par Note that the $L^2$ projection of $\bar{u}$ does not belong to $U_{ad}$. Instead, we will consider the standard interpolation operator $\piho: \Uad\longrightarrow \ddUad$.  
The following estimate (cf. \cite{FEM:BS2008}) is standard for $\bar{u}\in W^{1,\infty}(\Omega)$
\begin{equation}\label{eq:w1infesti}
    \|\bar{u}-\piho(\bar{u})\|^2_{L^2(T)}\le Ch^2_T\int_T |\nabla \bar{u}|^2\ \!dx\le Ch_T^2\|\bar{u}\|^2_{W^{1,\infty}(T)}\text{meas}(T)\le Ch_T^4\|\bar{u}\|^2_{W^{1,\infty}(T)}
\end{equation}
for any $T\in\cT_h$. However we need to modify the operator $\piho$ since it does not satisfy the orthogonal property with respect to the $L^2$ inner product as the operator $\Pi_h^0$ did in the previous section. 
When $h$ is sufficiently small, it is reasonable to assume there is no $T \in \mct$ such that $\min \limits_{\bar{T}} \olu = u_a$ and $\max \limits_{\bar{T}} \olu = u_b$ at the same time. We then define $\tluh \vert_T \in \ddUad$ as follows \cite{OCP:elliptic:rosch2005linear}:
\begin{align*}
    \tluh \vert_T:=
    \begin{cases}
        u_a,         & \min \limits_{\bar{T}} \olu = u_a \\
        u_b,         & \max \limits_{\bar{T}} \olu = u_b \;.\\
        \piho(\olu), & \text{otherwise}
    \end{cases}
\end{align*}

\begin{lemma}\label{eocp:sec-p1 approx control:lemma:auxiliary control satisfies continuous VI}
    For sufficiently small $h > 0$, we have 
    \begin{align}
    \big(\ol{p} + \beta \olu, u - \tluh \big)_{\ltO} \geq 0 \qquad \forall u \in U_{ad}.
    \end{align}
\end{lemma}
\begin{proof} 
    \par We decompose the domain $\Omega$ into three parts $\Ome = \Ome_a \cup \Ome_b \cup \mathcal{N}$, where 
    \begin{align*}
        \Ome_a &:= \{ x \in \Ome: \bar{u}(x) = u_a \}, \\
        \Ome_b &:= \{ x \in \Ome: \bar{u}(x) = u_b \}, \\
        \mathcal{N} &:= \Omega \setminus (\Ome_a \cup \Ome_b). 
    \end{align*}
    From \eqref{eocp:sec-optimality conditions:continuous KKT VI}, we conclude that $\ol{p} + \beta \olu \geq 0$ on $\Ome_a$, $\ol{p} + \beta \olu \leq 0$ on $\Ome_b$, and $\ol{p} + \beta \olu = 0$ on $\mathcal{N}$. 

    For any $u \in \Uad$ and $T \in \mathcal{T}_h$, we will show $\big(\ol{p} + \beta \olu, u - \tluh \big)_{T} \geq 0$. First, consider $T \in \mct$ in which there exists a $x \in T$ such that $\olu(x) = u_a$. Then, by definition, $\tluh = u_a $ on $T$. Thus, $u - \tluh \geq 0$ on $T$. For sufficiently small $h > 0$, we have $T \subset \Ome_a \cup \mathcal{N}$ and, thus, $\ol{p} + \beta \olu \geq 0$ on $T$. These imply $\big(\ol{p} + \beta \olu, u - \tluh \big)_{T} \geq 0$ on such a $T \in \mathcal{T}_h$. Similarly, consider $T \in \mct$ in which there exists a $x \in T$ such that $\olu(x) = u_b$. Then $u - \tluh = u - u_b \leq 0$ on $T$. In this case, we have $T \subset \Ome_b \cup \mathcal{N}$ and, thus, $\ol{p} + \beta \olu \leq 0$ on $T$ for sufficiently small $h > 0$. We also have $\big(\ol{p} + \beta \olu, u - \tluh \big)_{T} \geq 0$. Finally, consider $T \in \mct$ in which $\min \limits_{\bar{T}} \olu \neq u_a$ and $\max \limits_{\bar{T}} \olu \neq u_b$. Then $T \subset \mathcal{N}$. In this case, $\ol{p} + \beta \olu = 0$ and, thus, $\big(\ol{p} + \beta \olu, u - \tluh \big)_{T} = 0$. 
\end{proof} 

To obtain improved error estimates, we introduce the following sets:
        \begin{align*}
            \omct &:= \{T \in \mct: \olu = u_a \text{ or } \olu = u_b\}, \\
            \twmct &:= \{T \in \mct: u_a < \olu < u_b\}, \\
            \thmct &:= \mct \setminus \left(\omct \cup \twmct \right).
        \end{align*}
        \par Notice that for any $T \in \thmct$, there exists $x_1, x_2 \in T$ such that $\olu(x_1) = u_a$ or $\olu(x_1) = u_b$ and $u_a<\olu(x_2)<u_b.$ Additionally, we suppose that there exists $C$ independent of $h$ such that \cite{OCP:superconvergenceOfControl}
\begin{align}\label{assump:thmct}
    \text{meas}(\thmct) \leq C h.
\end{align}
\begin{theorem}\label{eocp:sec-p1 approx control:thm:L2 err in exact and discrete control}
    Let $\olu \in U_{ad} \cap W^{1,\infty}(\Omega)$ and $\ol{u}_h \in \ddUad \subset U_{ad}$ be the solutions of the problems \eqref{intro:sec-eocp:continuous minimization problem} and \eqref{eocp:sec-optimality conditions:discrete minimization}, respectively. Assume the assumption \eqref{assump:thmct} holds. Then there exists a constant C that depends on $|\ol{u}|_{W^{1,\infty}(\Omega)}$ and is independent of $h$ such that
    \begin{align}\label{higher order conv control}
        \|\ol{u} - \ol{u}_h \|_{L^2(\Omega)} \leq Ch^{\frac32} 
    \end{align}
    for sufficiently small $h > 0$.
\end{theorem}

\begin{proof}
    \par Using Lemma \ref{eocp:sec-p1 approx control:lemma:auxiliary control satisfies continuous VI} and choosing $u = \ol{u}_h \in \ddUad \subset U_{ad}$, we have
        \begin{align}
            0 &\leq \left(\ol{p}+\beta \ol{u},\ol{u}_h-\tluh\right)_{\ltO} \label{sec:p1 approx control:p1 proof continuous VI}.
        \end{align}
    By choosing $u_h = \tluh \in U^1_{ad,h}$ in the inequality \eqref{sec:err estimate:discrete VI operator form}, we have
    \begin{align}
        0 &\leq \big(-\left(\ol{p}_h + \beta \oluh\right) ,\ol{u}_h-\tluh \big)_{\ltO} \label{sec:p1 approx control:p1 proof discrete VI}.
    \end{align}
    \par Adding the inequalities \eqref{sec:p1 approx control:p1 proof continuous VI} and \eqref{sec:p1 approx control:p1 proof discrete VI}, we have
        \begin{align*}
              0 &\leq \big(\ol{p} - \ol{p}_h + \beta(\ol{u} - \ol{u}_h), \ol{u}_h - \tluh\big)_{\ltO} \notag \\
                  &= \left(\ol{p} - \ol{p}_h,\ol{u}_h - \tluh \right)_{\ltO} + \beta\left(\ol{u} - \ol{u}_h, \ol{u}_h - \tluh\right)_{\ltO} \notag \\
                  &=\left(\ol{p} - \ol{p}_h,\ol{u}_h - \tluh \right)_{\ltO} - \beta \|\ol{u} - \oluh\|^2_{\ltO} + \beta \left(\ol{u} - \oluh, \ol{u} - \tluh \right)_{\ltO}.
            \end{align*}
        \par Consequently,
            \begin{align}
                 \|\ol{u} - \oluh\|^2_{\ltO} &\leq\dfrac{1}{\beta}\left(\ol{p} - \ol{p}_h,\ol{u}_h - \tluh \right)_{\ltO} + \left(\ol{u} - \ol{u}_h, \ol{u}_h - \tluh\right)_{\ltO} \notag \\
                  &=\dfrac{1}{\beta}\left(\ol{p} - \ol{p}_h,\ol{u}_h - \ol{u} \right)_{\ltO} + \dfrac{1}{\beta}\left(\ol{p} - \ol{p}_h,\ol{u} - \tluh \right)_{\ltO} + \left(\ol{u} - \oluh, \ol{u} - \tluh \right)_{\ltO} \notag \\
                  &\leq \dfrac{1}{\beta}\left(\ol{p} - \ol{p}_h,\ol{u}_h - \ol{u} \right)_{\ltO} + \dfrac{1}{\beta}\left(\ol{p} - \ol{p}_h,\ol{u} - \tluh \right)_{\ltO} +  C\eps_1\|\ol{u} - \oluh\|^2_{\ltO} + \dfrac{C}{\eps_1}\|\ol{u} - \tluh \|^2_{\ltO} \notag \\
                  &=\underbrace{\dfrac{1}{\beta}\left(A^{\star}\left(A\ol{u} - y_d\right)-A_h^{\star}\big(A_h\ol{u}_h - y_d\right), \ol{u}_h-\ol{u}\big)_{\lt{\Ome}}}_{S_1} \notag \\
                  &\hspace{0.2in}+\underbrace{\dfrac{1}{\beta}\big(A^{\star}\left(A\ol{u} - y_d\right)-A_h^{\star}\left(A_h\ol{u}_h - y_d\right), \ol{u}-\tluh\big)_{\lt{\Ome}}}_{S_2} \notag\\
                  &\hspace{0.2in}+C\eps_1\|\ol{u} - \oluh\|^2_{\ltO} + \dfrac{C}{\eps_1}\|\ol{u} - \tluh \|^2_{\ltO}. \label{intermediatebound}
            \end{align}
        \par To estimate the terms $S_1$ and $S_2$, we follow the steps of the inequalities \eqref{boundT1}--\eqref{boundT2} to arrive at
        \begin{align}
            |S_1| &\leq Ch^4 (\|\ol{u}\|_{\ltO}^2 + \|y_d\|_{\ltO}^2) + \eps_2 \; \|\ol{u} - \oluh\|^2_{\ltO} \label{boundT7},\\
            |S_2| &\leq Ch^4 (\|\ol{u}\|_{\ltO}^2 + \|y_d\|_{\ltO}^2) + \eps_3 \; \|\ol{u} - \tluh\|^2_{\ltO}. \label{boundT8}
        \end{align}
        \par By combining the inequalities \eqref{intermediatebound}--\eqref{boundT8} and by appropriately choosing $\eps_1, \eps_2, \eps_3$, we have
        \begin{align} \label{link intermediate estimate}
            \|\ol{u} - \oluh\|^2_{\ltO} &\leq  Ch^4 (\|\ol{u}\|_{\ltO}^2 + \|y_d\|_{\ltO}^2) + C\|\ol{u} - \tluh\|^2_{\ltO}.
        \end{align}
        \par The remaining task is to estimate $\|\ol{u} - \tluh\|^2_{\ltO}.$ 
        Consider 
        \begin{align*}
            \|\olu - \tluh\|^2_{\ltO} = \underbrace{\sum_{T \in \omct} \|\olu - \tluh\|^2_{\lt{T}}}_{S_3} + \underbrace{\sum_{T \in \twmct} \|\olu - \tluh\|^2_{\lt{T}}}_{S_4} + \underbrace{\sum_{T \in \thmct} \|\olu - \tluh\|^2_{\lt{T}}}_{S_5}.
        \end{align*}
        \par By the definition of $\tluh$ and $\omct$, We have $S_3 = 0$. 
        \par To estimate $S_4$, it follows from the definition of $\tluh$ and \eqref{eq:w1infesti} that
        \[S_4 = \sum_{T \in \twmct} \|\olu - \tluh\|^2_{\lt{T}} =  \sum_{T \in \twmct} \|\olu - \piho(\olu)\|^2_{\lt{T}} \leq Ch^4 \|\olu\|^2_{W^{1,\infty}(\twmct)}.\]
        \par Before estimating $S_5$, for any $T \in \thmct$, consider $\hat{x}_T \in T$ such that, without loss of generality, $\olu(\hat{x}_T) = u_a$. Note that, for $T\in \cT_{h,3},$ there holds
        \begin{align*}
                    \|\ol{u} - \tluh\|^2_{\lt{T}}   &= \|\ol{u} - u_a\|^2_{\lt{T}} \\
                                                            &= \ds \int_{T} |\ol{u} - u_a|^2 \; dT \\
                                                            &= \ds \int_{T} |\ol{u}(x) - \ol{u}(\hat{x}_T)|^2 \; dT \\
                                                            &=  \ds \mathlarger{\mathlarger{\int}}_{T} \left|\sum_{|\alpha|=1} \left(x - \hat{x}_T\right)^\alpha \int_0^1\dfrac{1}{\alpha!}D^\alpha \olu\big(x + s(\hat{x}_T - x)\big)\;ds \right|^2 \; dT \\
                        &= \ds \int_{T} \left|\int_0^1 (x - \hat{x}_T) \cdot \nabla\olu\big(x + s(\hat{x}_T - x)\big) \; ds \right|^2 \; dT \\
                        &\leq \ds \int_{T} \int_0^1\left|(x - \hat{x_T})\right|^2 \; ds  \int_0^1\left|\nabla\olu\big(x + s(\hat{x}_T - x)\big)\right|^2 \; ds \; dT \\
                        &\leq C h_{T}^2 \ds  \int_0^1   \int_{T}\|\nabla\olu\|^2_{L^{\infty}(T)} \; dT \; ds \\
                        &\leq C h_{T}^2 \ds  \int_0^1   \|\nabla\olu\|^2_{L^{\infty}(T)} \int_{T} 1 \; dT \; ds \\
                        &\leq C h_{T}^2 \|\nab \olu \|^2_{L^{\infty}(T)} \text{meas}(T).
                    \end{align*}                   
        \par Using the assumption \eqref{assump:thmct}, we have 
        \begin{align}
            S_5 = \sum_{T \in \thmct} \|\olu - \tluh\|^2_{\lt{T}} 
            \leq Ch^{3} |\ol{u}|_{W^{1,\infty}(\thmct)}.
        \end{align}
        \par Ultimately, the desired estimate \eqref{higher order conv control} is obtained by combining the inequality \eqref{link intermediate estimate} with the bounds for $S_3,S_4,$ and $S_5.$
\end{proof}
\subsection{An a priori Error Estimate on the State and the Adjoint State}\label{eocp:sec-state}
\par This section is devoted to estimating $\energy{\ol{y} - \ol{y}_h}$ and $\energy{\ol{p} - \ol{p}_h}$. 
We will first establish an intermediate error estimate for $\energy{\ol{y} - y_h(\olu)}$, where $y_h(\olu)$ satisfies 
\begin{align}
a_h\big(y_h(\ol{u}), v_h\big) = (\ol{u},v_h)_{\cT_h} \qquad \forall v_h \in V_h.
\end{align} 

\begin{lemma}\label{sec:state:lemma:intermediate result}
    There exists a $C > 0$ independent of $h$ such that
    \begin{align*}
        \energy{ y_h(\olu) -  \ol{y}_h} \leq C \|\olu - \oluh \|_{\lt{\Ome}}.
    \end{align*}
\end{lemma}
\begin{proof}
    \par Note that, by definition,
    \begin{align*}
        \energy{ y_h(\olu) -  \ol{y}_h }^2
         & =
        a_h\big( y_h(\olu) -  \ol{y}_h,y_h(\olu) -  \ol{y}_h \big) \notag \\
         & =
        \big(\olu - \oluh,y_h(\olu) -  \ol{y}_h\big)_{\lt{\Ome}} \notag   \\
         & \leq
        \|\olu - \oluh \|_{\lt{\Ome}} \; \|y_h(\olu) -  \ol{y}_h\|_{\lt{\Ome}}.
    \end{align*}
    \par By Lemma \ref{sec:a priori analysis:lemma:DWDG Poincare},  we have
    \begin{align*}
        \|y_h(\olu) -  \ol{y}_h\|_{\lt{\Ome}} \leq C \energy{y_h(\olu) -  \ol{y}_h}.
    \end{align*}
    \par Hence, we have
    \begin{align}
        \energy{y_h(\olu) -  \ol{y}_h}^2 \leq C \|\olu - \oluh \|_{\lt{\Ome}} \; \energy{y_h(\olu) -  \ol{y}_h}. 
    \end{align}
    The desired result immediately follows.
\end{proof}
\begin{theorem} \label{sec:state:thm:energy err in exact and discrete state}
    Let $\ol{y} \in  H^2(\Ome)$ and $\ol{y}_h \in V_h$ be the solutions of the problems \eqref{intro:sec-eocp:continuous minimization problem} and \eqref{eocp:sec-optimality conditions:discrete minimization}, respectively. There exists a constant C that depends on $|\ol{y}|_{H^2(\Omega)}$ and is independent of $h$ such that \\
    \begin{align*}
        \energy{\ol{y} - \ol{y}_h} \leq Ch.
    \end{align*}
\end{theorem}
\begin{proof}
    \par Note that $\bar{y}_h(\bar{u})$ is the DWDG approximation of $\bar{y}$ defined by \eqref{eocp:sec-optimality conditions:continuous KKT PDE weakform}; hence, we have (cf. \cite{DWDG:LN2014}) $$\energy{\bar{y}-\bar{y}_h(\bar{u})}\le Ch.$$ Therefore, by the triangle inequality, we have
    \begin{equation}
        \energy{\bar{y}-\bar{y}_h}\le \energy{\bar{y}-\bar{y}_h(\bar{u})}+\energy{\bar{y}_h(\bar{u})-\bar{y}_h}\le Ch
    \end{equation}
    where we used Lemma \ref{sec:state:lemma:intermediate result}, Theorem \ref{sec:p0 approx control:thm:L2 err in exact and discrete control}, and Theorem \ref{eocp:sec-p1 approx control:thm:L2 err in exact and discrete control}.
\end{proof}
\par 
Next, we state a theorem that provides the error estimate for the adjoint state in the energy norm. First, we prove the following lemma:

\begin{lemma} \label{sec:adj:lemma:intermediate result}
    Let $\ol{y} \in  H^2(\Ome)$ be the solution of the problem \eqref{intro:sec-eocp:continuous minimization problem}. There exists a constant C that depends on $|\ol{y}|_{H^2(\Omega)}$ and is independent of $h$ such that \\
    \begin{align*}
        \energy{p_h(\ol{y}) - \ol{p}_h} \leq Ch,
    \end{align*}
    where  $p_h(\ol{y})$ satisfies $a_h\big(p_h(\ol{y}), v_h\big) = (y_{d} - \ol{y},v_h)_{\cT_h} \quad \text{for all } v_h \in V_h$.
\end{lemma}
\begin{proof}
    \par We have by the definition of $\energy{\cdot}$ that
    \begin{align*}
        \energy{p_h(\ol{y}) - \ol{p}_h}^2 & = a_h \big( p_h(\ol{y}) - \ol{p}_h , p_h(\ol{y}) - \ol{p}_h \big)                \\
                                          & = \big(y_{d} - \ol{y} - (y_{d} - \ol{y}_h), p_h(\ol{y}) - \ol{p}_h\big)_{L^2(\Ome)}  \\
                                          & = \big(\ol{y}_h - \ol{y}, p_h(\ol{y}) - \ol{p}_h\big)_{L^2(\Ome)}                        \\
                                          & \leq \|\ol{y}_h - \ol{y}\|_{L^2(\Ome)} \; \|p_h(\ol{y}) - \ol{p}_h\|_{L^2(\Ome)}.
    \end{align*}
    \par Notice that by Lemma \ref{sec:a priori analysis:lemma:DWDG Poincare}, $\|p_h(\ol{y}) - \ol{p}_h\|_{L^2(\Ome)} \leq \energy{p_h(\ol{y}) - \ol{p}_h}$.
    Subsequently, let $\Pi_h (\ol{y})$ be the standard nodal interpolant of $\ol{y}$ in $V_h$. By using the triangle inequality and Lemma \ref{sec:a priori analysis:lemma:DWDG Poincare}, we have
    \begin{align}
         \|\ol{y}_h - \ol{y}\|_{L^2(\Ome)} 
         & \leq \|\ol{y}_h -\Pi_h(\ol{y}) \|_{L^2(\Ome)} + \|\Pi_h(\ol{y}) - \ol{y}\|_{L^2(\Ome)} \notag                                                                                 \\
         & \leq \energy{\ol{y}_h -\Pi_h(\ol{y})}+ \|\Pi_h(\ol{y}) - \ol{y}\|_{L^2(\Ome)}    \notag                                                                                      \\
         & \leq\energy{\ol{y}_h - \ol{y}} + \energy{\ol{y} -\Pi_h(\ol{y})}+ \|\Pi_h(\ol{y}) - \ol{y}\|_{L^2(\Ome)}\label{sec:adj:consequence of DWDG Poincare and nodal interpolation}.
    \end{align}
    \par It is standard that $\|\Pi_h(\ol{y}) - \ol{y}\|_{L^2(\Ome)} \leq Ch^2 |\ol{y}|_{H^2(\Ome)}$ (cf. \cite{FEM:BS2008,FEM:ciarlet2002finite}) and $\energy{\ol{y} -\Pi_h(\ol{y})} \leq Ch |\ol{y}|_{H^2(\Ome)}$ (cf. \cite{EVI:DWDG:LRZ2020}).
     Finally, by Theorem \ref{sec:state:thm:energy err in exact and discrete state}, we have the desired estimate.
\end{proof}
\begin{theorem} \label{sec:adj:thm:energy err in exact and discrete adj}
    Let $\ol{p} \in  H^2(\Ome)$ and $\ol{p}_h \in V_h$ be the solutions to \eqref{eocp:sec-optimality conditions:continuous KKT adj weakform} and \eqref{eocp:sec-optimality conditions:discrete KKT adj weakform}, respectively. There exists a constant C that depends on $|\ol{p}|_{H^2(\Omega)}$ and is independent of $h$ such that \\
    \begin{align*}
        \energy{\ol{p} - \ol{p}_h} \leq Ch.
    \end{align*}
\end{theorem}
\begin{proof}
    Due to a similar argument as used in Theorem \ref{sec:state:thm:energy err in exact and discrete state} and the analysis presented in \cite{DWDG:LN2014, EVI:DWDG:LRZ2020}, it follows that $\energy{\ol{p} - p_h(\ol{y})} \leq Ch.$ The final result then follows from Lemma \ref{sec:adj:lemma:intermediate result} and the triangle inequality.
\end{proof}

\section{Numerical Experiments}\label{eocp:sec-numerical experiments}
\par In this section, we present several numerical examples to demonstrate the robustness of the proposed scheme and validate the theoretical results. These examples are generated using in-house MATLAB codes. The finite-dimensional problem obtained through the DWDG method is solved using the primal-dual active set algorithm \cite{PDAS:bergounioux1999primal,PDAS:hintermuller2002primal}.

\subsection{Example 1}\label{eocp:ex1}
\par In this example, we set $u_a = -\infty$, $u_b = \infty$, $y_d = (1+4\pi^4) \sin(\pi x_1) \sin(\pi x_2) $ and seek $(\ol{y}_h,\ol{u}_h,\ol{p}_h) \in V_h \times \Uadh \times V_h$ for $k \in \{0,1\}$ which solves the system \eqref{eocp:sec-optimality conditions:discrete KKT} on $\Ome = [0,1]^2$. This is an optimal control problem without control constraints.
\par The exact solution $(\ol{y},\ol{u},\ol{p}) \in H^2(\Ome) \times U_{ad} \times H^2(\Ome)$ is given by 
\begin{subequations}\label{sec:num results:ex1 exact solution}
    \begin{empheq}{align}
        \ol{y}(x_1,x_2) &= \sin(\pi x_1) \sin(\pi x_2) \label{sec:num results:ex1 exact y}\, ,
        \\
        \ol{u}(x_1,x_2) &= 2 \pi^2 \sin(\pi x_1) \sin(\pi x_2) \label{sec:num results:ex1 exact u}\, , \\
        \ol{p}(x_1,x_2) &= -2 \pi^2 \sin(\pi x_1) \sin(\pi x_2) \label{sec:num results:ex1 exact p}\, .
    \end{empheq}
\end{subequations}
\par Tables \ref{sec:num results:Table1:ex1_y_p0_u}, \ref{sec:num results:Table2:ex1_y_p1_u}, \ref{sec:num results:Table3:ex1_p_p0_u}, and \ref{sec:num results:Table4:ex1_p_p1_u} show the convergence of $\ol{y}_h$ to $\ol{y}$ and the convergence of $\ol{p}_h$ to $\ol{p}$ in the \emph{energy norm} for three different penalty parameters as $h \rightarrow 0$ respectively. We see a convergence of order 1 in the \emph{energy norm} as proved in Theorems \ref{sec:state:thm:energy err in exact and discrete state} and \ref{sec:adj:thm:energy err in exact and discrete adj}.

\begin{table}[!htb]
    \centering
    \begin{tabular}{|cc|cc|cc|cc|}
        \hline
        \multicolumn{2}{|c|}{}
                               & \multicolumn{2}{c|}{$\gamma = -1$}
                               & \multicolumn{2}{c|}{$\gamma = 0$}
                               & \multicolumn{2}{c|}{$\gamma = 5$}                                                        \\
        $h$                    & DOF                                & $\energy{\ol{y}_h - \ol{y}}$    & Rate & $\energy{\ol{y}_h - \ol{y}}$    & Rate & $\energy{\ol{y}_h - \ol{y}}$    & Rate \\
        \hline
        \multirow{1}{*}{1/2}   & 12                                 & 2.71e+00 & -- & 2.68e+00 & -- & 2.61e+00 & -- \\
        \multirow{1}{*}{1/4}   & 48                                 & 1.08e+00 & 1.33 & 1.10e+00 & 1.29 & 1.15e+00 & 1.18 \\
        \multirow{1}{*}{1/8}   & 192                                & 5.36e-01 & 1.01 & 5.46e-01 & 1.01 & 5.75e-01 & 1.00 \\
        \multirow{1}{*}{1/16}  & 768                                & 2.72e-01 & 0.98 & 2.77e-01 & 0.98 & 2.89e-01 & 0.99 \\
        \multirow{1}{*}{1/32}  & 3072                               & 1.38e-01 & 0.98 & 1.40e-01 & 0.99 & 1.45e-01 & 0.99 \\
        \multirow{1}{*}{1/64}  & 12288                              & 6.94e-02 & 0.99 & 7.03e-02 & 0.99 & 7.28e-02 & 1.00 \\
        \multirow{1}{*}{1/128} & 49152                              & 3.48e-02 & 0.99 & 3.53e-02 & 1.00 & 3.65e-02 & 1.00 \\
        \hline
    \end{tabular}
    \caption{Rates of convergence of $\energy{\ol{y}_h - \ol{y}}$ for Example \ref{eocp:ex1} using $\mathbb{P}_0$ approximation for $\ol{u}_h.$}
    \label{sec:num results:Table1:ex1_y_p0_u}
\end{table}

\begin{table}[!htb]
    \centering
    \begin{tabular}{|cc|cc|cc|cc|}
        \hline
        \multicolumn{2}{|c|}{}
                               & \multicolumn{2}{c|}{$\gamma = -1$}
                               & \multicolumn{2}{c|}{$\gamma = 0$}
                               & \multicolumn{2}{c|}{$\gamma = 5$}                                                        \\
        $h$                    & DOF                                & $\energy{\ol{y}_h - \ol{y}}$    & Rate & $\energy{\ol{y}_h - \ol{y}}$    & Rate & $\energy{\ol{y}_h - \ol{y}}$    & Rate \\
        \hline
        \multirow{1}{*}{1/2} & 12 & 6.24e-01 & --& 6.24e-01 & -- & 6.24e-01 & --  \\                                                                          
         \multirow{1}{*}{1/4} & 48 & 6.29e-01 & -0.01& 6.66e-01 & -0.09 & 7.69e-01 & -0.30  \\                                                                       
         \multirow{1}{*}{1/8} & 192 & 3.32e-01 & 0.93& 3.50e-01 & 0.93 & 3.95e-01 & 0.96  \\                                                                         
         \multirow{1}{*}{1/16} & 768 & 1.77e-01 & 0.91& 1.84e-01 & 0.93 & 2.02e-01 & 0.97  \\                                                                         
         \multirow{1}{*}{1/32} & 3072 & 9.16e-02 & 0.95& 9.47e-02 & 0.96 & 1.02e-01 & 0.98  \\                                                                        
         \multirow{1}{*}{1/64} & 12288 & 4.66e-02 & 0.97& 4.80e-02 & 0.98 & 5.15e-02 & 0.99  \\                                                                       
         \multirow{1}{*}{1/128} & 49152 & 2.35e-02 & 0.99& 2.42e-02 & 0.99 & 2.58e-02 & 1.00  \\ 
        \hline
    \end{tabular}
    \caption{Rates of convergence of $\energy{\ol{y}_h - \ol{y}}$ for Example \ref{eocp:ex1} using $\mathbb{P}_1$ approximation for $\ol{u}_h$.}
    \label{sec:num results:Table2:ex1_y_p1_u}
\end{table}

\begin{table}[!htb]
    \centering
    \begin{tabular}{|cc|cc|cc|cc|}
        \hline
        \multicolumn{2}{|c|}{}
                               & \multicolumn{2}{c|}{$\gamma = -1$}
                               & \multicolumn{2}{c|}{$\gamma = 0$}
                               & \multicolumn{2}{c|}{$\gamma = 5$}                                                           \\
        $h$                    & DOF                                & $\energy{\ol{p}_h - \ol{p}}$    & Rate  & $\energy{\ol{p}_h - \ol{p}}$    & Rate  & $\energy{\ol{p}_h - \ol{p}}$    & Rate  \\
        \hline
        \multirow{1}{*}{1/2}   & 12                                 & 1.07e+01 & --  & 1.07e+01 & --  & 1.07e+01 & --  \\
        \multirow{1}{*}{1/4}   & 48                                 & 1.16e+01 & -0.10 & 1.22e+01 & -0.19 & 1.41e+01 & -0.39 \\
        \multirow{1}{*}{1/8}   & 192                                & 6.45e+00 & 0.84  & 6.81e+00 & 0.84  & 7.67e+00 & 0.88  \\
        \multirow{1}{*}{1/16}  & 768                                & 3.48e+00 & 0.89  & 3.62e+00 & 0.91  & 3.98e+00 & 0.95  \\
        \multirow{1}{*}{1/32}  & 3072                               & 1.81e+00 & 0.94  & 1.87e+00 & 0.96  & 2.02e+00 & 0.98  \\
        \multirow{1}{*}{1/64}  & 12288                              & 9.20e-01 & 0.97  & 9.47e-01 & 0.98  & 1.02e+00 & 0.99  \\
        \multirow{1}{*}{1/128} & 49152                              & 4.64e-01 & 0.99  & 4.77e-01 & 0.99  & 5.10e-01 & 1.00  \\
        \hline
    \end{tabular}
    \caption{Rates of convergence of $\energy{\ol{p}_h - \ol{p}}$ for Example \ref{eocp:ex1} using $\mathbb{P}_0$ approximation for $\ol{u}_h$.}
    \label{sec:num results:Table3:ex1_p_p0_u}
\end{table}

\begin{table}[!htb]
    \centering
    \begin{tabular}{|cc|cc|cc|cc|}
        \hline
        \multicolumn{2}{|c|}{}
                               & \multicolumn{2}{c|}{$\gamma = -1$}
                               & \multicolumn{2}{c|}{$\gamma = 0$}
                               & \multicolumn{2}{c|}{$\gamma = 5$}                                                           \\
        $h$                    & DOF                                & $\energy{\ol{p}_h - \ol{p}}$    & Rate  & $\energy{\ol{p}_h - \ol{p}}$    & Rate  & $\energy{\ol{p}_h - \ol{p}}$    & Rate  \\
        \hline
        \multirow{1}{*}{1/2} & 12 & 1.07e+01 & --& 1.07e+01 & -- & 1.07e+01 & --  \\                                                                          
         \multirow{1}{*}{1/4} & 48 & 1.16e+01 & -0.10& 1.22e+01 & -0.19 & 1.41e+01 & -0.39  \\                                                                       
         \multirow{1}{*}{1/8} & 192 & 6.45e+00 & 0.84& 6.81e+00 & 0.84 & 7.67e+00 & 0.88  \\                                                                         
         \multirow{1}{*}{1/16} & 768 & 3.48e+00 & 0.89& 3.62e+00 & 0.91 & 3.98e+00 & 0.95  \\                                                                         
         \multirow{1}{*}{1/32} & 3072 & 1.81e+00 & 0.94& 1.87e+00 & 0.96 & 2.02e+00 & 0.98  \\                                                                        
         \multirow{1}{*}{1/64} & 12288 & 9.20e-01 & 0.97& 9.47e-01 & 0.98 & 1.02e+00 & 0.99  \\                                                                       
         \multirow{1}{*}{1/128} & 49152 & 4.64e-01 & 0.99& 4.77e-01 & 0.99 & 5.10e-01 & 1.00  \\
        \hline
    \end{tabular}
    \caption{Rates of convergence of $\energy{\ol{p}_h - \ol{p}}$ for Example \ref{eocp:ex1} using $\mathbb{P}_1$ approximation for $\ol{u}_h$.}
    \label{sec:num results:Table4:ex1_p_p1_u}
\end{table}
\par In Tables \ref{sec:num results:Table3:ex1_u p0} and \ref{sec:num results:Table4:ex1_u p1}, the convergence of $\ol{u}_h$ to $\ol{u}$ \eqref{sec:num results:ex1 exact u} as $h \rightarrow 0$ is shown. Following Theorem \ref{sec:p0 approx control:thm:L2 err in exact and discrete control}, we see an order 1 convergence for $\ol{u}_h \in \dUad$ in Table \ref{sec:num results:Table3:ex1_u p0}. A second-order convergence in the $L^2$ norm of $\ol{u}_h \in \ddUad$ is evident in Table \ref{sec:num results:Table4:ex1_u p1} across the three distinct penalty parameters. This result is expected since $\ol{u} \in H^2(\Ome)$ as there are no constraints on the control variable.

\begin{table}[!htb]
    \centering
    \begin{tabular}{|cc|cc|cc|cc|}
        \hline
        \multicolumn{2}{|c|}{}
                               & \multicolumn{2}{c|}{$\gamma = -1$}
                               & \multicolumn{2}{c|}{$\gamma = 0$}
                               & \multicolumn{2}{c|}{$\gamma = 5$}                                                        \\
        $h$                    & DOF                                & $\|\ol{u}_h - \ol{u}\|_{\ltO}$    & Rate & $\|\ol{u}_h - \ol{u}\|_{\ltO}$    & Rate & $\|\ol{u}_h - \ol{u}\|_{\ltO}$    & Rate \\
        \hline
        \multirow{1}{*}{1/2}   & 04                                 & 3.26e+00 & -- & 3.26e+00 & -- & 3.26e+00 & -- \\
        \multirow{1}{*}{1/4}   & 16                                 & 2.67e+00 & 0.29 & 2.74e+00 & 0.25 & 2.91e+00 & 0.17 \\
        \multirow{1}{*}{1/8}   & 64                                 & 1.31e+00 & 1.03 & 1.32e+00 & 1.05 & 1.34e+00 & 1.11 \\
        \multirow{1}{*}{1/16}  & 256                                & 6.49e-01 & 1.01 & 6.50e-01 & 1.02 & 6.53e-01 & 1.04 \\
        \multirow{1}{*}{1/32}  & 1024                               & 3.23e-01 & 1.00 & 3.24e-01 & 1.01 & 3.24e-01 & 1.01 \\
        \multirow{1}{*}{1/64}  & 4096                               & 1.62e-01 & 1.00 & 1.62e-01 & 1.00 & 1.62e-01 & 1.00 \\
        \multirow{1}{*}{1/128} & 16384                              & 8.08e-02 & 1.00 & 8.08e-02 & 1.00 & 8.08e-02 & 1.00 \\
        \hline
    \end{tabular}
    \caption{Rates of convergence of $\|\ol{u}_h - \ol{u}\|_{\ltO}$ for Example \ref{eocp:ex1} using $\mathbb{P}_0$ approximation for $\ol{u}_h$.}
    \label{sec:num results:Table3:ex1_u p0}
\end{table}

\begin{table}[!htb]
    \centering
    \begin{tabular}{|cc|cc|cc|cc|}
        \hline
        \multicolumn{2}{|c|}{}
                               & \multicolumn{2}{c|}{$\gamma = -1$}
                               & \multicolumn{2}{c|}{$\gamma = 0$}
                               & \multicolumn{2}{c|}{$\gamma = 5$}                                                        \\
        $h$                    & DOF                                & $\|\ol{u}_h - \ol{u}\|_{\ltO}$    & Rate & $\|\ol{u}_h - \ol{u}\|_{\ltO}$    & Rate & $\|\ol{u}_h - \ol{u}\|_{\ltO}$    & Rate \\
        \hline
        \multirow{1}{*}{1/2}   & 12                                 & 2.29e+00 & -- & 2.29e+00 & -- & 2.29e+00 & -- \\
        \multirow{1}{*}{1/4}   & 48                                 & 1.20e+00 & 0.93 & 1.22e+00 & 0.91 & 1.35e+00 & 0.76 \\
        \multirow{1}{*}{1/8}   & 192                                & 2.70e-01 & 2.15 & 2.84e-01 & 2.10 & 3.33e-01 & 2.02 \\
        \multirow{1}{*}{1/16}  & 768                                & 6.80e-02 & 1.99 & 7.18e-02 & 1.98 & 8.42e-02 & 1.98 \\
        \multirow{1}{*}{1/32}  & 3072                               & 1.75e-02 & 1.96 & 1.84e-02 & 1.96 & 2.14e-02 & 1.98 \\
        \multirow{1}{*}{1/64}  & 12288                              & 4.47e-03 & 1.97 & 4.70e-03 & 1.97 & 5.40e-03 & 1.99 \\
        \multirow{1}{*}{1/128} & 49152                              & 1.13e-03 & 1.98 & 1.19e-03 & 1.99 & 1.36e-03 & 1.99 \\
        \hline
    \end{tabular}
    \caption{Rates of convergence of $\|\ol{u}_h - \ol{u}\|_{\ltO}$ for Example \ref{eocp:ex1} using $\mathbb{P}_1$ approximation for $\ol{u}_h$.}
    \label{sec:num results:Table4:ex1_u p1}
\end{table}
\subsection{Example 2}\label{eocp:ex2}
\par In the second example, we still seek $(\ol{y}_h,\ol{u}_h,\ol{p}_h) \in V_h \times \Uadh \times V_h$ for $k \in \{0,1\}$ solving the system \eqref{eocp:sec-optimality conditions:discrete KKT} on $\Ome = [0,1]^2$, where $y_d = (1+4\pi^4) \sin(\pi x_1) \sin(\pi x_2)$. However, we set $u_a = 3$ and $u_b = 15$ as in \cite{OCP:elliptic:rosch2005linear}. Note that the box constraints on the control variable, in this case, are non-trivial.
\par The exact solution $(\ol{y},\ol{u},\ol{p}) \in H^2(\Ome) \times U_{ad} \times H^2(\Ome)$ is given by 
\begin{subequations}\label{sec:num results:ex2 exact solution}
    \begin{empheq}{align}
        \ol{y}(x_1,x_2) &= \sin(\pi x_1) \sin(\pi x_2) \label{sec:num results:ex2 exact y} \, ,
        \\
        \ol{u}(x_1,x_2) &=
        \begin{cases}
            u_a,                                & \text{if } 2 \pi^2 \sin(\pi x_1) \sin(\pi x_2) < u_a,          \\
            2 \pi^2 \sin(\pi x_1) \sin(\pi x_2), & \text{if } 2 \pi^2 \sin(\pi x_1) \sin(\pi x_2) \in [u_a, u_b], \\
            u_b,                                & \text{if } 2 \pi^2 \sin(\pi x_1) \sin(\pi x_2) > u_b \, ,
        \end{cases} \label{sec:num results:ex2 exact u}\\
        \ol{p}(x_1,x_2) &= -2 \pi^2 \sin(\pi x_1) \sin(\pi x_2) \label{sec:num results:ex2 exact p} \, .
    \end{empheq}
\end{subequations}
\par
\begin{figure}[H]
    \centering
    \begin{minipage}{.365\textwidth}
        \hspace{-0.7in}
        \includegraphics[width=1.30\linewidth]{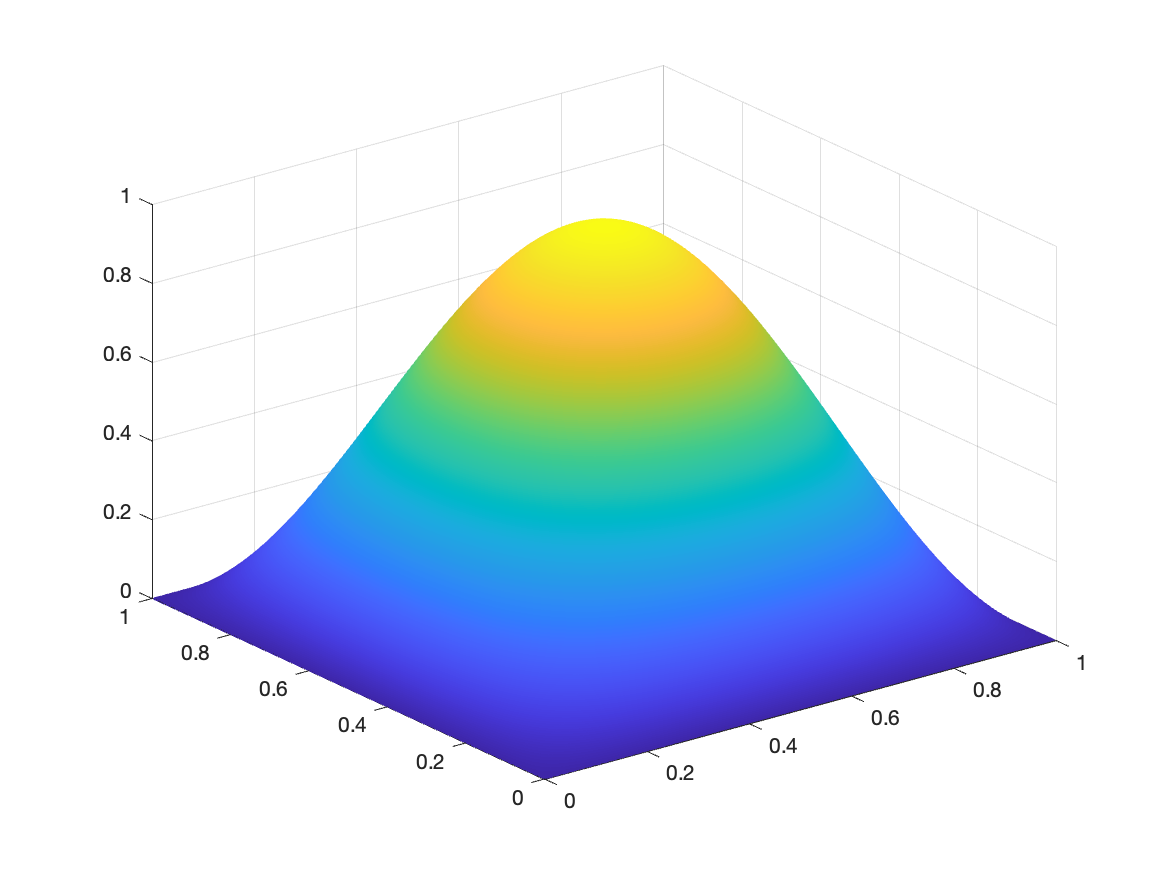}
    \end{minipage}%
    \begin{minipage}{.365\textwidth}
        \centering
        \includegraphics[width=1.30\linewidth]{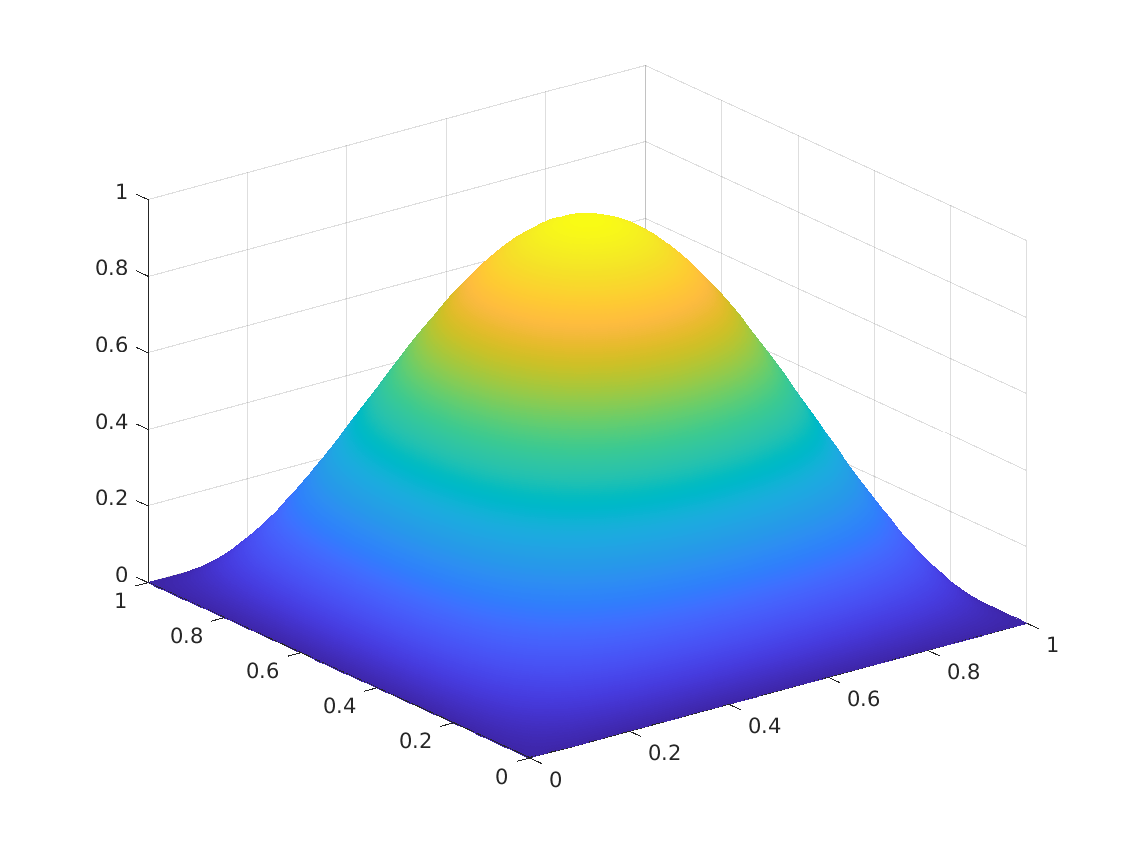}
    \end{minipage}
    \caption{Results for Example \ref{eocp:ex2}: $\ol{y}_h$ (left) and $\ol{y}$ (right); $h = \frac{1}{128}$.}
    \label{fig:eocp:ex2-y}
\end{figure}

\begin{figure}[H]
    \centering
    \begin{minipage}{.365\textwidth}
        \hspace{-0.7in}
        \includegraphics[width=1.30\linewidth]{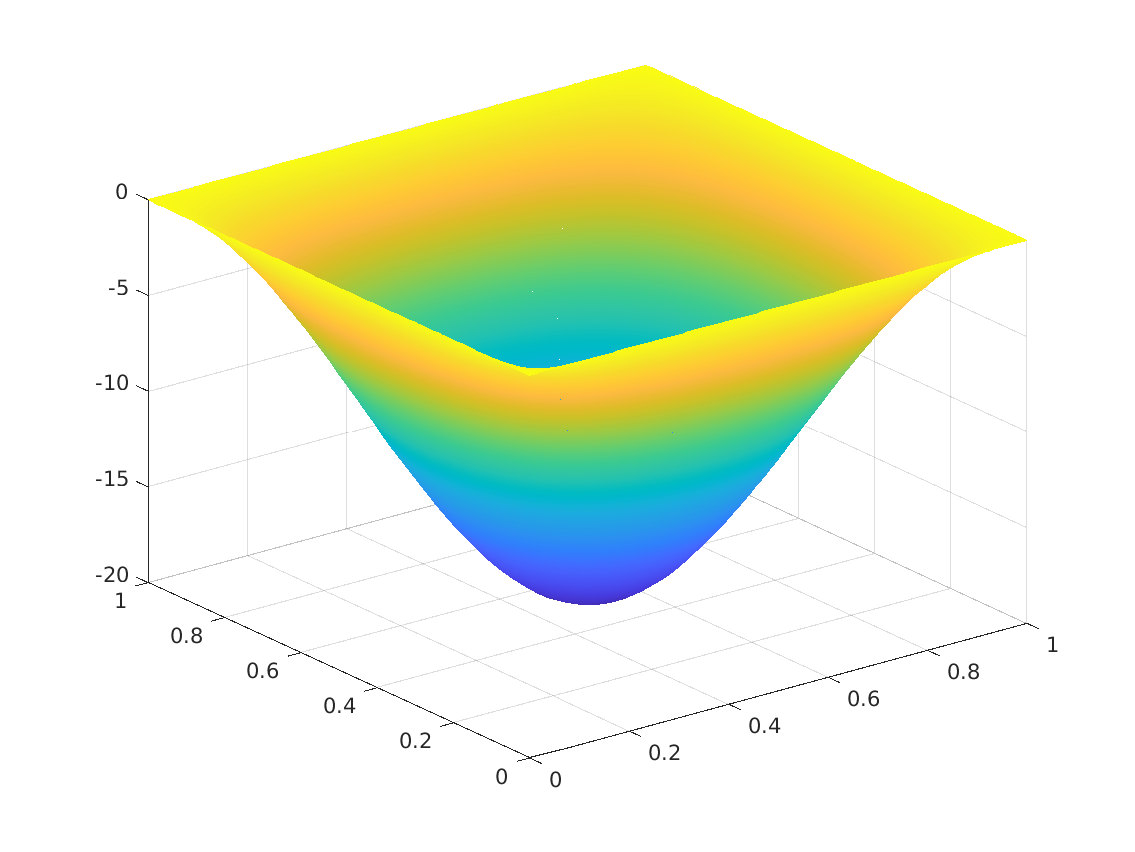}
    \end{minipage}%
    \begin{minipage}{.365\textwidth}
        \centering
        \includegraphics[width=1.30\linewidth]{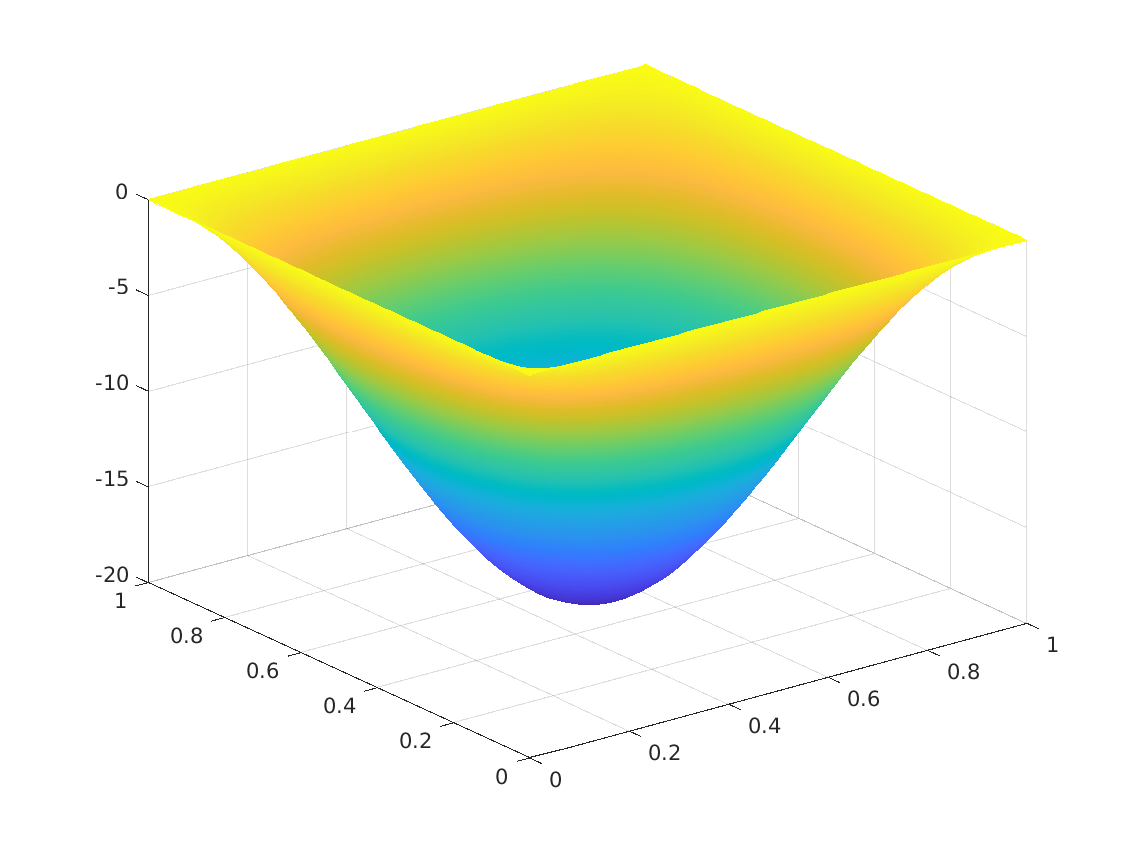}
    \end{minipage}
    \caption{Results for Example \ref{eocp:ex2}: $\ol{p}_h$ (left) and $\ol{p}$ (right); $h = \frac{1}{128}$.}
    \label{fig:eocp:ex2-p}
\end{figure}

\begin{figure}[H]
    \centering
    \begin{minipage}{.365\textwidth}
        \hspace{-0.7in}
        \includegraphics[width=1.30\linewidth]{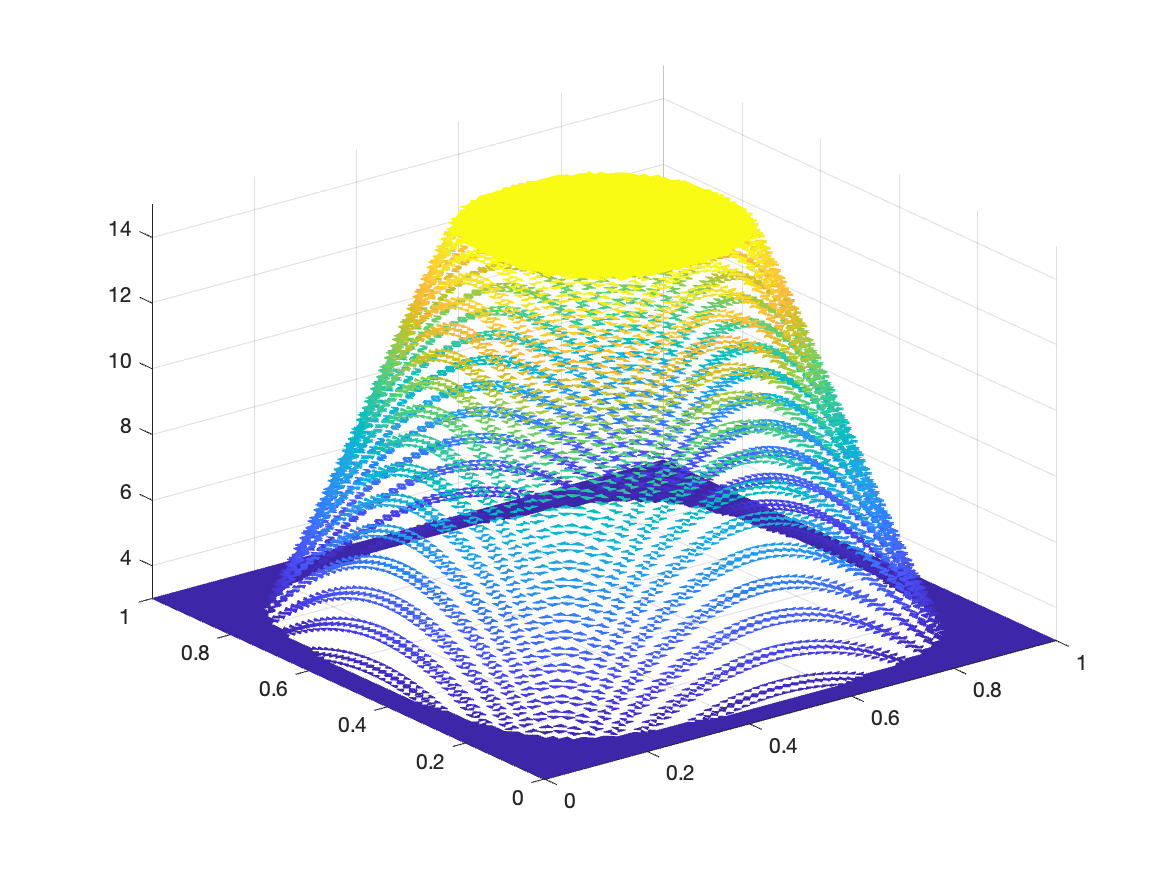}
    \end{minipage}%
    \begin{minipage}{.365\textwidth}
        \centering
        \includegraphics[width=1.30\linewidth]{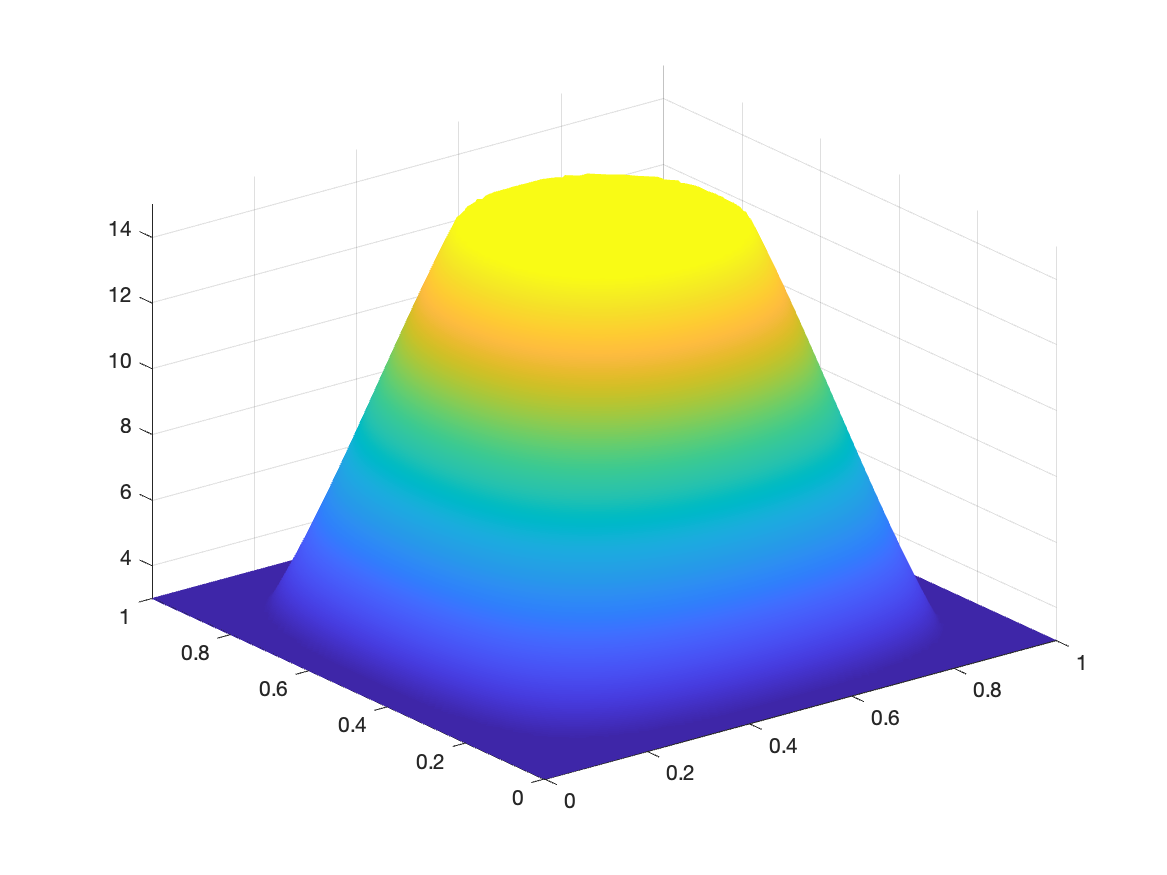}
    \end{minipage}
    \caption{Results for Example \ref{eocp:ex2}: $\ol{u}_h \in \dUad$ (left), $\ol{u}_h \in \ddUad$ (right), $h = \frac{1}{128}$.}
    \label{fig:eocp:ex2-u} 
\end{figure}

 As in the case of Example \ref{eocp:ex1}, the errors in the energy norm and the orders of convergence of $\ol{y}_h$ to $\ol{y}$ and $\ol{p}_h$ to $\ol{p}$ as $h \rightarrow 0$ are tabulated in Tables \ref{sec:num results:Table5:ex2_y_p0_u}, \ref{sec:num results:Table6:ex2_y_p1_u}, \ref{sec:num results:Table7:ex2_p_p0_u}, and \ref{sec:num results:Table8:ex2_p_p1_u} respectively. Also, in Tables \ref{sec:num results:Table7:ex2_u p0} and \ref{sec:num results:Table8:ex2_u p1}, the orders of convergence of $\ol{u}_h \in \dUad$ and $\ol{u}_h \in \ddUad$ to $\ol{u}$  in the $L^2$ norm are shown. The numerical findings in this example align with the theoretical results established in Section \ref{sec:a priori analysis}. 
Furthermore, we plot the discrete and continuous optimal state and adjoint state in Figure~\ref{fig:eocp:ex2-y} and Figure~\ref{fig:eocp:ex2-p}. The discrete $\mathbb{P}_0$ and $\mathbb{P}_1$ approximation of the optimal control are also shown in Figure~\ref{fig:eocp:ex2-u}. 

\begin{table}[!htb]
    \centering
    \begin{tabular}{|cc|cc|cc|cc|}
        \hline
        \multicolumn{2}{|c|}{}
                               & \multicolumn{2}{c|}{$\gamma = -1$}
                               & \multicolumn{2}{c|}{$\gamma = 0$}
                               & \multicolumn{2}{c|}{$\gamma = 5$}                                                        \\
        $h$                    & DOF                                & $\energy{\ol{y}_h - \ol{y}}$    & Rate & $\energy{\ol{y}_h - \ol{y}}$    & Rate & $\energy{\ol{y}_h - \ol{y}}$    & Rate \\
        \hline
        \multirow{1}{*}{1/2}   & 12                                 & 2.89e+00 & -- & 2.86e+00 & -- & 2.78e+00 & -- \\
        \multirow{1}{*}{1/4}   & 48                                 & 1.02e+00 & 1.50 & 1.04e+00 & 1.46 & 1.08e+00 & 1.36 \\
        \multirow{1}{*}{1/8}   & 192                                & 5.20e-01 & 0.97 & 5.30e-01 & 0.97 & 5.58e-01 & 0.96 \\
        \multirow{1}{*}{1/16}  & 768                                & 2.63e-01 & 0.98 & 2.68e-01 & 0.98 & 2.80e-01 & 0.99 \\
        \multirow{1}{*}{1/32}  & 3072                               & 1.34e-01 & 0.98 & 1.36e-01 & 0.98 & 1.41e-01 & 0.99 \\
        \multirow{1}{*}{1/64}  & 12288                              & 6.73e-02 & 0.99 & 6.82e-02 & 0.99 & 7.08e-02 & 1.00 \\
        \multirow{1}{*}{1/128} & 49152                              & 3.38e-02 & 0.99 & 3.42e-02 & 1.00 & 3.54e-02 & 1.00 \\
        \hline
    \end{tabular}
    \caption{Rates of convergence of $\energy{\ol{y}_h - \ol{y}}$ for Example \ref{eocp:ex2} using $\mathbb{P}_0$ approximation for $\ol{u}_h.$}
    \label{sec:num results:Table5:ex2_y_p0_u}
\end{table}

\begin{table}[!htb]
    \centering
    \begin{tabular}{|cc|cc|cc|cc|}
        \hline
        \multicolumn{2}{|c|}{}
                               & \multicolumn{2}{c|}{$\gamma = -1$}
                               & \multicolumn{2}{c|}{$\gamma = 0$}
                               & \multicolumn{2}{c|}{$\gamma = 5$}                                                        \\
        $h$                    & DOF                                & $\energy{\ol{y}_h - \ol{y}}$    & Rate & $\energy{\ol{y}_h - \ol{y}}$    & Rate & $\energy{\ol{y}_h - \ol{y}}$    & Rate \\
        \hline
        \multirow{1}{*}{1/2} & 12 & 5.87e-01 & --& 5.91e-01 & -- & 6.02e-01 & --  \\                                                                          
         \multirow{1}{*}{1/4} & 48 & 6.22e-01 & -0.08& 6.60e-01 & -0.16 & 7.65e-01 & -0.35  \\                                                                       
         \multirow{1}{*}{1/8} & 192 & 3.28e-01 & 0.92& 3.46e-01 & 0.93 & 3.90e-01 & 0.97  \\                                                                         
         \multirow{1}{*}{1/16} & 768 & 1.76e-01 & 0.89& 1.84e-01 & 0.91 & 2.02e-01 & 0.95  \\                                                                         
         \multirow{1}{*}{1/32} & 3072 & 9.16e-02 & 0.95& 9.46e-02 & 0.96 & 1.02e-01 & 0.98  \\                                                                        
        \multirow{1}{*}{1/64} & 12288 & 4.66e-02 & 0.97& 4.80e-02 & 0.98 & 5.15e-02 & 0.99  \\                                                                       
         \multirow{1}{*}{1/128} & 49152 & 2.35e-02 & 0.99& 2.42e-02 & 0.99 & 2.58e-02 & 1.00  \\    
        \hline
    \end{tabular}
    \caption{Rates of convergence of $\energy{\ol{y}_h - \ol{y}}$ for Example \ref{eocp:ex2} using $\mathbb{P}_1$ approximation for $\ol{u}_h.$}
    \label{sec:num results:Table6:ex2_y_p1_u}
\end{table}

\begin{table}[!htb]
    \centering
    \begin{tabular}{|cc|cc|cc|cc|}
        \hline
        \multicolumn{2}{|c|}{}
                               & \multicolumn{2}{c|}{$\gamma = -1$}
                               & \multicolumn{2}{c|}{$\gamma = 0$}
                               & \multicolumn{2}{c|}{$\gamma = 5$}                                                           \\
        $h$                    & DOF                                & $\energy{\ol{p}_h - \ol{p}}$    & Rate  & $\energy{\ol{p}_h - \ol{p}}$    & Rate  & $\energy{\ol{p}_h - \ol{p}}$    & Rate  \\
        \hline
        \multirow{1}{*}{1/2}   & 12                                 & 1.07e+01 & --  & 1.07e+01 & --  & 1.07e+01 & --  \\
        \multirow{1}{*}{1/4}   & 48                                 & 1.16e+01 & -0.11 & 1.22e+01 & -0.19 & 1.41e+01 & -0.39 \\
        \multirow{1}{*}{1/8}   & 192                                & 6.45e+00 & 0.84  & 6.81e+00 & 0.84  & 7.67e+00 & 0.88  \\
        \multirow{1}{*}{1/16}  & 768                                & 3.48e+00 & 0.89  & 3.62e+00 & 0.91  & 3.98e+00 & 0.95  \\
        \multirow{1}{*}{1/32}  & 3072                               & 1.81e+00 & 0.94  & 1.87e+00 & 0.96  & 2.02e+00 & 0.98  \\
        \multirow{1}{*}{1/64}  & 12288                              & 9.20e-01 & 0.97  & 9.47e-01 & 0.98  & 1.02e+00 & 0.99  \\
        \multirow{1}{*}{1/128} & 49152                              & 4.64e-01 & 0.99  & 4.77e-01 & 0.99  & 5.10e-01 & 1.00  \\
        \hline
    \end{tabular}
    \caption{Rates of convergence of $\energy{\ol{p}_h - \ol{p}}$ for Example \ref{eocp:ex2} using $\mathbb{P}_0$ approximation for $\ol{u}_h$ .}
    \label{sec:num results:Table7:ex2_p_p0_u}
\end{table}

\begin{table}[!htb]
    \centering
    \begin{tabular}{|cc|cc|cc|cc|}
        \hline
        \multicolumn{2}{|c|}{}
                               & \multicolumn{2}{c|}{$\gamma = -1$}
                               & \multicolumn{2}{c|}{$\gamma = 0$}
                               & \multicolumn{2}{c|}{$\gamma = 5$}                                                           \\
        $h$                    & DOF                                & $\energy{\ol{p}_h - \ol{p}}$    & Rate  & $\energy{\ol{p}_h - \ol{p}}$    & Rate  & $\energy{\ol{p}_h - \ol{p}}$    & Rate  \\
        \hline
        \multirow{1}{*}{5.00e-01} & 12 & 1.07e+01 & --& 1.07e+01 & -- & 1.07e+01 & --  \\                                                                          
         \multirow{1}{*}{2.50e-01} & 48 & 1.16e+01 & -0.10& 1.22e+01 & -0.19 & 1.41e+01 & -0.39  \\                                                                       
         \multirow{1}{*}{1.25e-01} & 192 & 6.45e+00 & 0.84& 6.81e+00 & 0.84 & 7.67e+00 & 0.88  \\                                                                         
         \multirow{1}{*}{6.25e-02} & 768 & 3.48e+00 & 0.89& 3.62e+00 & 0.91 & 3.98e+00 & 0.95  \\                                                                         
         \multirow{1}{*}{3.12e-02} & 3072 & 1.81e+00 & 0.94& 1.87e+00 & 0.96 & 2.02e+00 & 0.98  \\                                                                        
         \multirow{1}{*}{1.56e-02} & 12288 & 9.20e-01 & 0.97& 9.47e-01 & 0.98 & 1.02e+00 & 0.99  \\                                                                       
         \multirow{1}{*}{7.81e-03} & 49152 & 4.64e-01 & 0.99& 4.77e-01 & 0.99 & 5.10e-01 & 1.00  \\
        \hline
    \end{tabular}
    \caption{Rates of convergence of $\energy{\ol{p}_h - \ol{p}}$ for Example \ref{eocp:ex2} using $\mathbb{P}_1$ approximation for $\ol{u}_h$ .}
    \label{sec:num results:Table8:ex2_p_p1_u}
\end{table}

\begin{table}[htb]
    \centering
    \begin{tabular}{|cc|cc|cc|cc|}
        \hline
        \multicolumn{2}{|c|}{}
                               & \multicolumn{2}{c|}{$\gamma = -1$}
                               & \multicolumn{2}{c|}{$\gamma = 0$}
                               & \multicolumn{2}{c|}{$\gamma = 5$}                                                        \\
        $h$                    & DOF                                & $\|\ol{u}_h - \ol{u}\|_{\ltO}$    & Rate & $\|\ol{u}_h - \ol{u}\|_{\ltO}$    & Rate & $\|\ol{u}_h - \ol{u}\|_{\ltO}$    & Rate \\
        \hline
        \multirow{1}{*}{1/2}   & 04                                 & 3.26e+00 & -- & 3.26e+00 & -- & 3.26e+00 & -- \\
        \multirow{1}{*}{1/4}   & 16                                 & 2.67e+00 & 0.29 & 2.74e+00 & 0.25 & 2.91e+00 & 0.17 \\
        \multirow{1}{*}{1/8}   & 64                                 & 1.09e+00 & 1.30 & 1.10e+00 & 1.32 & 1.11e+00 & 1.38 \\
        \multirow{1}{*}{1/16}  & 256                                & 5.33e-01 & 1.03 & 5.34e-01 & 1.04 & 5.37e-01 & 1.05 \\
        \multirow{1}{*}{1/32}  & 1024                               & 2.62e-01 & 1.02 & 2.62e-01 & 1.03 & 2.63e-01 & 1.03 \\
        \multirow{1}{*}{1/64}  & 4096                               & 1.36e-01 & 0.95 & 1.36e-01 & 0.95 & 1.36e-01 & 0.95 \\
        \multirow{1}{*}{1/128} & 16384                              & 6.82e-02 & 1.00 & 6.82e-02 & 1.00 & 6.82e-02 & 1.00 \\
        \hline
    \end{tabular}
    \caption{Rates of convergence of $\|\ol{u}_h - \ol{u}\|_{\ltO}$ for Example \ref{eocp:ex2} using $\mathbb{P}_0$ approximation for $\ol{u}_h.$}
    \label{sec:num results:Table7:ex2_u p0}
\end{table}

\begin{table}[!htb]
    \centering
    \begin{tabular}{|cc|cc|cc|cc|}
        \hline
        \multicolumn{2}{|c|}{}
                               & \multicolumn{2}{c|}{$\gamma = -1$}
                               & \multicolumn{2}{c|}{$\gamma = 0$}
                               & \multicolumn{2}{c|}{$\gamma = 5$}                                                        \\
        $h$                    & DOF                                & $\|\ol{u}_h - \ol{u}\|_{\ltO}$    & Rate & $\|\ol{u}_h - \ol{u}\|_{\ltO}$    & Rate & $\|\ol{u}_h - \ol{u}\|_{\ltO}$    & Rate \\
        \hline
        \multirow{1}{*}{1/2}   & 12                                 & 2.57e+00 & -- & 2.57e+00 & -- & 2.57e+00 & -- \\
        \multirow{1}{*}{1/4}   & 48                                 & 1.16e+00 & 1.14 & 1.18e+00 & 1.12 & 1.28e+00 & 1.01 \\
        \multirow{1}{*}{1/8}   & 192                                & 5.96e-01 & 0.96 & 5.99e-01 & 0.98 & 6.10e-01 & 1.07 \\
        \multirow{1}{*}{1/16}  & 768                                & 1.73e-01 & 1.79 & 1.73e-01 & 1.79 & 1.75e-01 & 1.80 \\
        \multirow{1}{*}{1/32}  & 3072                               & 4.43e-02 & 1.96 & 4.45e-02 & 1.96 & 4.50e-02 & 1.96 \\
        \multirow{1}{*}{1/64}  & 12288                              & 2.20e-02 & 1.01 & 2.20e-02 & 1.01 & 2.21e-02 & 1.02 \\
        \multirow{1}{*}{1/128} & 49152                              & 7.92e-03 & 1.47 & 7.92e-03 & 1.47 & 7.94e-03 & 1.48 \\
        \hline
    \end{tabular}
    \caption{Rates of convergence of $\|\ol{u}_h - \ol{u}\|_{\ltO}$ for Example \ref{eocp:ex2} using $\mathbb{P}_1$ approximation for $\ol{u}_h.$}
    \label{sec:num results:Table8:ex2_u p1}
\end{table}
\section{Summary}\label{sec:Conclusions} 
\par In this work, we studied a dual-wind Discontinuous Galerkin (DWDG) scheme to discretize an optimal control problem constrained by box constraints on the control, with the governing PDE represented by Poisson's equation. This discretization process led to a finite-dimensional optimization problem, which was solved using the primal-dual active set algorithm. We established the error estimates \emph{a priori} in the appropriate norms for the solution ($\ol{y}$, $\ol{u}$, $\ol{p}$).

Several numerical tests were conducted to demonstrate error convergence in suitable norms. Potential future research is to improve the convergence rate of the discrete control variable to the exact control in the \(L^2\) norm by making use of a projection operator and the discrete adjoint variable \(p_h\) in a post-processing step, following the approach outlined in \cite{OCP:superconvergenceOfControl}.

\par Furthermore, we plan to extend this research by developing a new DG method based on the DG finite element differential calculus \cite{DWDG:FLN2016}, for when the PDE constraint is a convection-diffusion equation within a convection-dominated regime (cf. \cite{liu2024multigrid,heinkenschloss2010local,leykekhman2012local}). This will allow us to establish refined \emph{a priori} error estimates for such problems. It is also interesting to consider fast solvers for DWDG (cf. \cite{brenner2020multigrid,liu2025robust}) and DWDG for optimal control problems with pointwise state constraints (cf. \cite{brenner2021p1,liu2024discontinuous,brenner2023multigrid}).

\bigskip
\noindent
\textbf{Acknowledgements} 
    The first author gratefully acknowledges the Department of Applied Mathematics at Florida Polytechnic University for providing funds that contributed to the successful completion of this research.
    The work of the second author was partially supported 
    by the National Science Foundation under Grant 
    DMS-2111059.  
    The work of the fourth author was partially supported 
    by the National Science Foundation under Grant 
    DMS-2111004.

\section*{References}

\bibliographystyle{plain}
\bibliography{references}

\end{document}